\documentclass[preprint, review, 12pt]{elsarticle}

\usepackage{amssymb, amsmath, bm}

\journal{Journal of Non-Newtonian Fluid Mechanics}

\usepackage{enumerate}
\usepackage[dvipsnames]{xcolor}

\newcommand{\De}{\text{De}}
\newcommand{\Wi}{\text{Wi}}

\begin{document}

\begin{frontmatter}



\title{The Method of Harmonic Balance for the Giesekus Model under Oscillatory Shear}


\author[iitk]{Shivangi Mittal}
\author[iitk]{Yogesh M. Joshi}
\author[fsu]{Sachin Shanbhag}

\affiliation[iitk]{organization={Department of Chemical Engineering, Indian Institute of Technology},
            addressline={}, 
            city={Kanpur},
            postcode={208016}, 
            state={Uttar Pradesh},
            country={India}}

\affiliation[fsu]{organization={Department of Scientific Computing, Florida State University},
            addressline={Dirac Science Library}, 
            city={Tallahassee},
            postcode={32306}, 
            state={Florida},
            country={United States}}

\begin{abstract}
The method of harmonic balance (HB) is a spectrally accurate method used to obtain periodic steady state solutions to dynamical systems subjected to periodic perturbations. We adapt HB to solve for the stress response of the Giesekus model under  large amplitude oscillatory shear (LAOS) deformation. HB transforms the system of differential equations to a set of nonlinear algebraic equations in the Fourier coefficients. Convergence studies find that the difference between the HB and true solutions decays exponentially with the number of harmonics ($H$) included in the ansatz as $e^{-m H}$. The decay coefficient $m$ decreases with increasing strain amplitude, and exhibits a ``U" shaped dependence on applied frequency. The computational cost of HB increases slightly faster than linearly with $H$. The net result of rapid convergence and modest increase in computational cost with increasing $H$ implies that HB outperforms the conventional method of using numerical integration to solve differential constitutive equations under oscillatory shear. Numerical experiments find that HB is simultaneously about three orders of magnitude cheaper, and several orders of magnitude more accurate  than numerical integration. Thus, it offers a compelling value proposition for parameter estimation or model selection.
\end{abstract}


\begin{highlights}
\item HB transforms a system of ordinary differential equations describing a dynamical system subjected to periodic perturbations into a set of nonlinear algebraic equations in frequency space.
\item For the Giesekus model under large amplitude oscillatory shear, we find that harmonic balance is orders of magnitude cheaper and more accurate than the conventional method of using numerical integration to solve differential constitutive equations.
\end{highlights}

\begin{keyword}
LAOS \sep spectral method \sep Fourier series \sep numerical method



\end{keyword}

\end{frontmatter}


\section{Introduction}

Rheological studies aim to understand the physics underlying a material's flow behavior, and quantify it in terms of measurable
flow parameters. Stress relaxation experiments, creep tests, and oscillatory shear measurements form the backbone of this enterprise. Among these, oscillatory shear experiments allow us to probe systems over time scales that span decades without complications that arise from steep ramps and abrupt step inputs. Here, a periodic shear strain $\gamma(t)=\gamma_{0}\sin\left(\omega t\right)$ with amplitude $\gamma_{0}$ and angular frequency $\omega$ is imposed, and the stresses generated in the material are measured.

When $\gamma_{0}$ is small, the resulting small amplitude oscillatory shear (SAOS) tests do not disturb the equilibrium structure, and characterize the linear response of the material. As $\gamma_{0}$ and $\omega$ increase, the behavior becomes increasingly nonlinear. Such large amplitude oscillatory shear (LAOS) measurements have been extensively used to study rheological phenomena including shear thinning/thickening and strain softening/hardening \cite{ma2020effects, khandavalli2015large, chan2018nonlinear, wapperom2005numerical}, time-dependent structural buildup or breakage \cite{li2009nonlinearity, dimitriou2013describing, min2014microstructure, armstrong2016dynamic, armstrong2021simple, Donley2019}, pseudoplasticity and elastoviscoplasticity \cite{ewoldt2010large,  stickel2013response, armstrong2016dynamic}, shear banding \cite{dimitriou2012rheo, goudoulas2017nonlinearities, kate2012large, radhakrishnan2018shear}, wall slip \cite{atalik2004occurrence,yang2017dynamic, klein2007separation, graham1995wall}, gelation \cite{suman2022large, kim2014microstructure, sun2015large, ng2011large}, chain stretch and entanglement in polymeric systems \cite{wagner2011analysis, hyun2007fourier, hoyle2014large, cho2010scaling, cho2015effect, cho2016viscoelasticity}, etc. Several analytical approaches have been introduced to interpret experimental LAOS data including Fourier series \cite{macdonald1969rheological}, power series \cite{pearson1982behavior}, Pade approximants \cite{giacomin2015pade}, Chebyshev polynomials \cite{Ewoldt2008}, stress decomposition methods \cite{cho2005geometrical, bae2017analytical}, characterstic waveforms \cite{klein2007separation}, weakly nonlinear intrinsic parameters \cite{hyun2009establishing},  sequence of physical processes \cite{rogers2012sequence}, etc.

\subsection{Representations for Oscillatory Shear Experiments}

In SAOS, stress profiles becomes sinusoidal after initial transients decay. In this periodic steady (or alternance) state, the shear stress becomes \cite{ferry1980viscoelastic} ,
\begin{equation}
\sigma_{12}^{\text{SAOS}}(t)=\gamma_{0}\left(G'(\omega)\sin\omega t+G''(\omega)\cos\omega t\right)\label{eq:SAOS}
\end{equation}
where $G'(\omega)$ and $G''(\omega)$ are the storage and loss moduli, respectively. The corresponding normal stress differences,
$N_{1}=\sigma_{11}-\sigma_{22}$ and $N_{2}=\sigma_{22}-\sigma_{33}$, mimic the upper convected Maxwell
(UCM) model \cite{ferry1980viscoelastic},
\begin{equation}
\begin{split}
\dfrac{N_{1}^{\text{SAOS}}(t)}{\gamma_{0}^{2}}  = G'(\omega) + & \left(G''(\omega)-\dfrac{G''(2\omega)}{2}\right)\sin2\omega t +  \\
& \left(-G'(\omega)+\dfrac{G''(2\omega)}{2}\right)\cos2\omega t,\\
\end{split}
\label{eq:N1_ucm}
\end{equation}
and $N_{2}^{\text{SAOS}}(t)  =0$. As $\gamma_{0}$ and $\omega$ are increased, the steady state stresses remain periodic, but are no longer perfect sinusoids
as higher harmonics are excited. Using symmetry, this nonlinear response can
be represented by a Fourier series with only odd harmonics for the
shear stress \cite{giacomin1993large, nam2008prediction},
\begin{equation}
\sigma_{12}(t)=\gamma_{0}\sum_{\substack{n=1\\
\text{odd}
}
}^{\infty}G_{n}'(\omega,\gamma_{0})\sin n\omega t+G_{n}''(\omega,\gamma_{0})\cos n\omega t,\label{eq:shear}
\end{equation}
where $G_{n}^{\prime}$ and $G_{n}^{\prime\prime}$ are sine and cosine Fourier coefficients, respectively, corresponding to the $n$th harmonic of the shear stress. Similarly, the normal stress differences contain only even harmonics \cite{cho2016viscoelasticity},
\begin{equation}
N_{1}(t)=\gamma_{0}^{2} \sum_{\substack{n=0\\\text{even}}}^{\infty} F_{n}'(\omega,\gamma_{0})\sin n\omega t+F_{n}''(\omega,\gamma_{0})\cos n\omega t,\label{eq:firstNormal}
\end{equation}
\begin{equation}
N_{2}(t)=\gamma_{0}^{2}\sum_{\substack{n=0\\\text{even}}}^{\infty} S_{n}'(\omega,\gamma_{0})\sin n\omega t+S_{n}''(\omega,\gamma_{0})\cos n\omega t,\label{eq:secondNormal}
\end{equation}
where $F_{n}^{\prime}$ and $F_{n}^{\prime\prime}$ ($S_{n}^{\prime}$
and $S_{n}^{\prime\prime}$) are sine and cosine Fourier coefficients,
respectively, corresponding to the $n$th harmonic of the first (second)
normal stress difference.

In the power series representation, a further Taylor series
expansion in $\gamma_{0}$ is performed. In the medium amplitude oscillatory shear (MAOS) or asymptotically nonlinear regime, the Taylor series is truncated after the $\gamma_{0}^{3}$ term. In this regime \cite{pearson1982behavior, giacomin1993large, cho2016viscoelasticity},
\begin{align}
\sigma_{12}(t) & =\gamma_{0}\left(G_{11}'\sin\omega t+G_{11}''\cos\omega t\right)+\nonumber \\
 & \gamma_{0}^{3}\left(G_{31}'\sin\omega t+G_{31}''\cos\omega t+G_{33}'\sin3\omega t+G_{33}''\cos3\omega t\right)+\mathcal{O}\left(\gamma_{0}^{5}\right)\label{eq:shear_power_series}
\end{align}
The power series coefficients $G_{ij}^{\prime}(\omega)$ and $G_{ij}^{\prime\prime}(\omega)$ are functions of only frequency, and are often called intrinsic coefficients. The coefficients $G_{11}^{\prime}\equiv G^{\prime}$ and $G_{11}^{\prime\prime}\equiv G^{\prime\prime}$ are identical to the linear viscoelastic moduli. The first index of the subscript ($i$) corresponds to the power of $\gamma_{0}$, while the second index ($j$) corresponds to the harmonic. The Fourier coefficients
are related to these power series coefficients via $G_{j}^{\prime}=\sum_{i=j}^{\infty} \gamma_{0}^{i-1}G_{ij}^{\prime}$
and $G_{j}^{\prime\prime}=\sum_{i=j}^{\infty} \gamma_{0}^{i-1}G_{ij}^{\prime\prime}$. Similarly, the power series representations of $N_{1}(t)$ and $N_{2}(t)$ in the MAOS regime are \cite{cho2016viscoelasticity},
\begin{equation}
N_{1}(t)=F_{00}''+\gamma_{0}^{2}\left(F_{20}''+F_{22}'\sin2\omega t+F_{22}''\cos2\omega t\right)+\mathcal{O}\left(\gamma_{0}^{4}\right)\label{eq:n1_power_series}
\end{equation}
\begin{equation}
N_{2}(t)=S_{00}''+\gamma_{0}^{2}\left(S_{20}''+S_{22}'\sin2\omega t+S_{22}''\cos2\omega t\right)+\mathcal{O}\left(\gamma_{0}^{4}\right)\label{eq:n2_power_series}
\end{equation}
where the intrinsic coefficients $F_{ij}^{\prime}$ ($S_{ij}^{\prime}$) and $F_{ij}^{\prime\prime}$ ($S_{ij}^{\prime\prime}$) are functions of only frequency, and are analogous to $G_{ij}^{\prime}$ and $G_{ij}^{\prime\prime}$, respectively. A table of the symbols used in this work is provided in supplementary material.

MAOS moduli have been widely studied theoretically and experimentally. Similar to the linear viscoelastic $G'$ and $G''$, the MAOS moduli $G_{33}'$ and $G_{33}''$ also follow nonlinear Kramers-Kronig relations \cite{shanbhag2022kramers1}, which can be used to assess the validity of experimental MAOS data \cite{shanbhag2022kramers2}. In the MAOS regime, the strain amplitude is barely large enough to trigger the second term of the power series but not other higher order terms. The ratio of third to first harmonic intensity is given by: 
\begin{equation}
I_{3/1} = \frac{\sqrt{G_{3}'^{2} + G_{3}''^{2}}}{\sqrt{G_{1}'^{2} + G_{1}''^{2}}} \approx \gamma_{0}^{2}\frac{\sqrt{(G_{33}')^{2}+(G_{33}'')^{2}}}{\sqrt{(G_{11}')^{2}+(G_{11}'')^{2}}}.
\end{equation}
The above approximation is valid only in the MAOS regime (small $\gamma_{0}$) since the second term in $G_{1}^{\prime} =G_{11}^{\prime} + \gamma_{0}^{2}G_{31}^{\prime}$ vanishes as $\gamma_{0}\rightarrow0$. This leads to a quadratic dependence given by $I_{3/1}\propto\gamma_{0}^{2}$.

\subsection{Constitutive Models}

Constitutive models are often necessary to attribute physical
meaning and interpret experimental rheological data. They consider the molecular forces generated within a material upon deformation, and yield differential or integral equations that relate stress and strain history. Numerous studies have sought to find the oscillatory shear response of different constitutive theories. Thus far, exact analytical expressions for shear and normal stresses corresponding to oscillatory strain inputs of arbitrary strain amplitude have been obtained for only two nontrivial constitutive equations: the corotational Maxwell model \cite{saengow2015exact, poungthong2019exact}, and the Oldroyd 8-constant model \cite{saengow2017exact, saengow2017normal}. However, both these models are quasilinear; the constitutive equations do not contain any nonlinear terms in stress.

Most analytical solutions for nonlinear constitutive equations have been obtained for only the MAOS regime, starting with early work on the Doi-Edwards theory for entangled polymer melts \cite{helfand1982calculation}. Indeed the list of constitutive equations with analytical expressions for the MAOS regime includes the Giesekus model \cite{nam2008prediction, kate2012large}, molecular stress function theory \cite{wagner2011analysis}, rigid rod-like polymer model \cite{bharadwaj2015constitutive, bird2014dilute}, Curtiss-Bird model \cite{fan1984kinetic}, fourth order fluids \cite{bharadwaj2014general}, simple emulsions \cite{yu2002modeling}, time-strain separable Kaye-Bernstein-Kerarsley-Zapas model \cite{martinetti2019time}, etc. By design, these analytical expressions do not provide the full LAOS response. They are applicable only for a narrow range of $\gamma_{0}$, in which $\gamma_0$ is large enough to elicit the weakest nonlinear modes, but small enough not to excite higher harmonics. In the absence of the full LAOS solution, it becomes difficult to apply these results to real systems.

The lack of analytical solutions at arbitrary strain amplitudes leaves numerical solution of the differential or integral constitutive equations as the only alternative. However, the performance of these methods degrades significantly at large $\gamma_{0}$ and $\omega$. Unlike typical numerical quadrature or time-stepping algorithms, spectral methods avoid these pitfalls. For example, we recently proposed a spectral method for time-strain separable integral equations \cite{shanbhag2021spectral}, which is not only more accurate than numerical quadrature, but also two to three orders of magnitude faster. It derives its speed and accuracy from a Fourier series representation of a portion of the integral equation. This combination of speed and accuracy facilitates its use in applications like Bayesian inference of model parameters \cite{suman2022large}, which requires the model to be evaluated thousands of times \cite{Shanbhag2010, Takeh2011}. The goal of this work is to adopt a similar technique for differential constitutive models in LAOS.

In this work, we demonstrate this approach for the Giesekus model \cite{giesekus1982simple}, which is a popular Maxwell-type differential
constitutive equation. It was originally developed for polymer solutions
\cite{yoo1989steady, schleiniger1991remark, yao1998extensional}, where
the extra stress tensor $\bm{\sigma}$, 
\begin{equation}
\bm{\sigma}=\bm{\sigma}_{p}+\bm{\sigma}_{s},\label{eq:totalStress}
\end{equation}
comprises contributions from the polymer ($\bm{\sigma}_{p}$), and 
Newtonian solvent, $\bm{\sigma}_{s} = \eta_{s}\dot{\bm{\gamma}}$, where
$\dot{\bm{\gamma}}$ is the deformation tensor and $\eta_{s}$ is
solvent viscosity. Subsequently, it has been applied to wormlike micellar surfactant solutions \cite{holz1999shear, fischer1997non, rehage2015experimental, bandyopadhyay2005effect, kate2012large}, protein dispersions \cite{kokini2000integral, Dhanasekharan2001}, concentrated dispersions \cite{duvarci2017saos}, etc. The polymer contribution is governed by, 
\begin{equation}
\frac{1}{\lambda}\bm{\sigma}_{p}+\overset{\nabla}{\bm{\sigma}}_{p}+\frac{\alpha}{\lambda G}\bm{\sigma}_{p} \cdot \bm{\sigma}_{p} = G \dot{\bm{\gamma}},\label{eq:tensor_form}
\end{equation}
where $\lambda$ and $G$ are the relaxation time and modulus, and
$\alpha$ is the anisotropy parameter ($0<\alpha<1$) that controls nonlinearity. $\overset{\nabla}{\bm{\sigma}}_{p}$ is the upper convected derivative defined as: 
\begin{equation}
\overset{\nabla}{\bm{\sigma}}_{p}=\frac{\partial\bm{\sigma}_{p}}{\partial t}+ \bm{v} \cdot \nabla\bm{\sigma}_{p} - \left(\bm{\nabla v}\right)^{T} \cdot \bm{\sigma}_{p}-\bm{\sigma}_{p} \cdot \bm{\nabla v},
\end{equation}
where $\bm{v}$ is fluid velocity. For $\alpha=0$, the Giesekus model reduces to the UCM model with a solvent contribution, and is therefore similar to the Jeffreys model. Coupled and uncoupled multimode extensions of the Giesekus model have been successfully used for predicting complex rheological phenomena for many systems in different flow fields \cite{calin2010determination, quinzani1990modeling, debbaut2002large, oztekin1994quantitative, atalik2002non, borzacchiello2016orientation}.

Early studies on this model revealed its capability to capture
shear thinning in a steady shear flow \cite{yoo1989steady}. Holtz
et al. analytically solved the response of the Giesekus model subjected to step shear strain of magnitude $\gamma$, and verified its use for a surfactant solution \cite{holz1999shear}. The stress relaxation modulus $G(t,\gamma)=\sigma_{12}(t)/\gamma$ was found to be, 
\begin{equation}
G(t,\gamma)=\frac{Ge^{-t/\lambda}}{1+2\alpha^{2}\gamma^{2}\left[1-\cosh\left(t/\lambda\right)\right]e^{-t/\lambda}+\alpha\gamma^{2}\left(1-e^{-t/\lambda}\right)}.
\end{equation}
Although the Giesekus model is not time-strain separable in general, it becomes separable in the long time limit ($t\rightarrow\infty$),
$G(t,\gamma)=G(t)h(\gamma)$ with $h(\gamma)=(1+\alpha(1-\alpha)\gamma^{2})^{-1}$.
The weakly nonlinear response of the Giesekus model has also been
investigated. Nam et al. \cite{nam2008prediction} obtained analytical expressions for normal stress differences truncated beyond $\gamma_{0}^{2}$.
Ewoldt et al. proposed a Chebyshev polynomial-based method to describe the intra-cycle behavior of Giesekus model fingerprints \cite{Ewoldt2008}. Ewoldt and McKinley further established a mathematical and physical framework to understand self-intersections in viscous Lissajous curves for the Giesekus model under LAOS \cite{ewoldt2010secondary}. Their analysis attributed this feature to the presence of strong elastic nonlinearities that emerge in LAOS. Calin et al.  proposed a method to obtain the parameters of a multimode Giesekus model from experimental SAOS and LAOS data for a semi-dilute polymer solution, using the third harmonic to obtain the anisotropy parameter $\alpha$ \cite{calin2010determination}. Bae and Cho proposed a semi-analytical approach for understanding MAOS response of the Giesekus model \cite{bae2015semianalytical}. Rogers and Lettinga explained the LAOS response of the Giesekus model as a sequence of physical processes, a unique approach that is quite different from traditional Fourier analysis \cite{rogers2012sequence}.

 The linear viscoelastic response of the Giesekus model is identical to the Maxwell model; thus $G_{11}'/G=\De^{2}/(1+\De^{2})$ and $G_{11}''/G=\De/(1+\De^{2})$, where the Deborah number $\De = \omega \lambda$. Gurnon and Wagner studied the effect of the anisotropy parameter on the stress response in LAOS and further confirmed the presence of secondary loops in Lissajous curves even for a wormlike micellar solution \cite{kate2012large}. They also reported the MAOS solution of the Giesekus model in terms of intrinsic power series coefficients. The intrinsic third harmonic terms for shear stress are given by, 
\begin{align}
\dfrac{G_{31}'}{G} & =\dfrac{\alpha\De^{4}\left(-21-41\De^{2}-8\De^{4}+4\alpha\left(4+7\De^{2}\right)\right)}{4\left(1+\De^{2}\right)^{3}\left(1+4\De^{2}\right)}\label{eq:G1p_gurnon}\\
\dfrac{G_{31}''}{G} & =-\dfrac{\alpha\De^{3}\left(9+11\De^{2}-10\De^{4}+2\alpha\left(-3-\De^{2}+8\De^{4}\right)\right)}{4\left(1+\De^{2}\right)^{3}\left(1+4\De^{2}\right)}\label{eq:G11pp_gurnon}\\
\dfrac{G_{33}'}{G} & =\dfrac{\alpha\De^{4}\left(-21+30\De+51\De^{4}+4\alpha\left(4-17\De^{2}+3\De^{4}\right)\right)}{4\left(1+\De^{2}\right)^{3}\left(1+4\De^{2}\right)\left(1+9\De^{2}\right)}\\
\dfrac{G_{33}''}{G} & =\dfrac{\alpha\De^{3}\left(-3+48\De^{2}+33\De^{4}-18\De^{6}+\alpha\left(2-48\De^{2}+46\De^{4}\right)\right)}{4\left(1+\De^{2}\right)^{3}\left(1+4\De^{2}\right)\left(1+9\De^{2}\right)}.
\end{align}
They also obtained similar expressions for the intrinsic power series coefficients of the normal stresses $\sigma_{11}$ and $\sigma_{22}$ corresponding to the zeroth and second harmonics \cite{kate2012large}.

\subsection{Initial Value Problem and Harmonic Balance}

The Giesekus model, like a majority of models for complex
fluids, is expressed as a differential constitutive equation (see
eqn \ref{eq:tensor_form}). It describes a set of partial differential
equations (PDEs) that can be plugged into a momentum balance equation
to compute stress and flow fields for arbitrary deformations. However,
in oscillatory shear experiments, the induced strain field is ideally homogeneous across the sample. In this scenario, the PDEs reduce to a set of nonlinear ordinary differential equations (ODEs), which can be posed as an initial value problem (IVP).

Typically, these nonlinear IVPs are solved numerically using
a time-stepping algorithm like Runge-Kutta \cite{Atkinson2009}. After initial transients decay, the solution eventually reaches a periodic limit cycle or alternance state, which is often the primary object of interest. This approach has several advantages: 

\begin{enumerate}[(i)]
\item \textbf{generality}: it can be easily adapted for
arbitrary constitutive models, 
\item \textbf{software}: sophisticated algorithms for solving
IVPs are already available, 
\item \textbf{evolution}: it tracks the evolution of the
system from the initial state to the alternance state, similar to
experiments. 
\end{enumerate}
However, it also comes with some disadvantages, especially
if only the final periodic solution is of interest: 
\begin{enumerate}[(i)]
\item \textbf{accuracy}: stiffness of these equations especially
at large $\gamma_{0}$ and $\omega$ can pose numerical problems,
which can be addressed by using a fine time step or accepting the numerical damping that comes with a coarse step. Implicit solvers are useful in dealing with stiff ODEs, but come with additional computational cost.  
\item \textbf{computational cost}: in addition to the increased
cost that arises from demands for greater accuracy, initial transients
can take a while to decay; we cannot algorithmically speed up this
process. 
\end{enumerate}
On the other hand, harmonic balance (HB) is a powerful technique
to efficiently solve periodic dynamical systems. Such systems are
encountered across various domains of science and engineering such
as transport and energy systems, acoustics, weather patterns, mechanical
dampers, alternating current powered processes \cite{krack2019harmonic},  etc.
When such systems are subjected to oscillatory inputs, they eventually
reach a periodic steady state. This periodic response can be expressed
as a Fourier series, i.e. a linear combination of sinusoidal waveforms.
This Fourier series representation or ansatz is plugged into the system of nonlinear differential equations describing the dynamics. We can determine the unknown coefficients of the Fourier series by matching or \textit{balancing}
them for each \textit{harmonic} up to a desired level of truncation; this step gives HB its name.

Interestingly, the \textit{numerical} technique of HB is closely related to previous efforts to \textit{analytically} determine the MAOS response of different constitutive models. However, extending the analytical approach to higher harmonics required to describe the LAOS response is tedious, and results in unwieldy expressions. Thus, HB can be thought of as a numerical generalization of this analytical approach. In this paper, we demonstrate how HB can be used to accurately and efficiently obtain the LAOS response of the Giesekus model from
first principles, by exploiting appropriate symmetries.

\subsection{Layout and Scope}

In this work, we establish the method of HB for solving differential constitutive relations subjected to LAOS in a computationally efficient and spectrally accurate manner by overcoming the shortcomings of conventional analytical and numerical techniques. We introduce the HB framework in section \ref{sec:methods}, and specialize it for solving the Giesekus model in LAOS. In section \ref{sec:results}, we demonstrate the viability of the computational protocol, and study the convergence properties of HB. Comparison with the conventional approach of numerical integration (NI) indicates that in terms of both speed and accuracy, HB is superior by orders of magnitude.

\section{Methods}
\label{sec:methods}

We begin by introducing the HB framework for a general system of ODEs. We then adapt it for the Giesekus model by exploiting symmetry constraints specific to oscillatory rheology. We review the conventional approach of NI for solving the IVP, and outline a protocol for comparing results with HB.

\subsection{Harmonic Balance for Constitutive Modeling}
\label{subsec:Harmonic-Balance}

The simplest way to introduce HB is to consider an algebraic differential equation in a single dependent variable $y(t)$, 
\begin{equation}
f(y,\dot{y},t)=0,\label{eq:intro_eqn_system}
\end{equation}
where $t$ is the independent variable, and $\dot{y}=dy/dt$. For
a periodic output, $y(t)=y(t+T)$, where $T$ is the period of oscillation. The Fourier series representation
of $y(t)$ up to $H$ harmonics is, 
\begin{equation}
y(t)\approx y_{H}(t)=\sum_{k=-H}^{H}Y_{k}\,e^{ik\omega t},\label{eq:FourierSeries}
\end{equation}
where frequency $\omega=2\pi/T$. Here $Y_{k}$ denotes the complex
Fourier coefficient associated with the $k$th harmonic, which is
defined as 
\begin{equation}
Y_{k}=\dfrac{1}{T}\int_{0}^{T}y(t)\,e^{-ik\omega t}dt.\label{eq:FourierTransform}
\end{equation}
We can express the summation in equation \eqref{eq:FourierSeries} as, 
\begin{equation}
y_{H}(t)=\bm{M}\bm{h}(\omega t),\label{eq:yh}
\end{equation}
where $\bm{M}=[Y_{-H},\cdots,Y_{0},\cdots,Y_{H}]$ is a row vector
of $2H+1$ Fourier coefficients, and $\bm{h}(\omega t)$ is the corresponding
column vector of Fourier basis functions, 
\begin{equation}
\bm{h}(\omega t)=\left[e^{-iH\omega t},\cdots,e^{-i\omega t},1,e^{i\omega t},\cdots,e^{iH\omega t}\right]^{T}\label{eq:h}
\end{equation}
Eqn (\ref{eq:yh}) effectively describes a dot product, due to the
shapes of $\bm{M}$ and $\bm{h}(\omega t)$. However, thinking of
it as a matrix product is helpful with generalization to systems of
equations later. Differentiating equation \eqref{eq:FourierSeries}, 
\begin{equation}
\dot{y}_{H}(t)=\sum_{k=-H}^{H}(ik\omega)Y_{k}\,e^{ik\omega t},\label{eq:FourierSeriesDvt}
\end{equation}
we find that the Fourier coefficient of the derivative $\dot{y}_{H}(t)$
corresponding to $k$th harmonic is $(ik\omega)Y_{k}$. Thus, we can
arrange these coefficients into a row vector with $2H+1$ elements
$\dot{\bm{M}}=\bm{MD}$, where $\bm{D}=\text{diag}(-i\omega H,\cdots,0,\cdots,i\omega H)$
is a $(2H+1)\times(2H+1)$ diagonal matrix. This relation is instrumental in transforming differential equations in the time domain into algebraic equations in the frequency domain using Fourier transforms. When the Fourier ansatz is substituted into equation (\ref{eq:intro_eqn_system}), we get: 
\begin{equation}
f(y_{H},\dot{y}_{H},t)=r_{H}(t),\label{eq:f_and_residual}
\end{equation}
where we typically obtain a nonzero residual $r_{H}$, due to truncation
of Fourier series after $H$ harmonics. $r_{H}(t)$ inherits periodicity
from $y_{H}(t)$, and can also be described using a Fourier series
truncated after $H$ harmonics. It can be represented by the matrix
multiplication, 
\begin{equation}
r_{H}(t)=\bm{R}\bm{h}(\omega t),\label{eq:residual}
\end{equation}
where $\bm{R}$ is a row vector of Fourier coefficients of $r_{H}(t)$
with $2H+1$ elements. Our original goal was to solve eqn \eqref{eq:intro_eqn_system}. In HB, this is approximated by requiring $\bm{R}=\bm{0}$. From equations (\ref{eq:residual}) and \eqref{eq:f_and_residual}, this is equivalent to solving $f(y_{H},\dot{y}_{H},t) = 0$. Note that $\bm{R}$ is a function of the Fourier coefficients $\bm{M}$. Practically, setting $\bm{R}(\bm{M})=\bm{0}$ involves solving a nonlinear algebraic system of $2H+1$ unknown Fourier coefficients that are the elements of $\bm{M}$.

Generalization of HB to a system of $N$ differential
equations $\bm{f}(\bm{y},\bm{\dot{y}},t)=\bm{0}$ is conceptually
straightforward. Here, $\bm{y}(t)$ and $\bm{y}_{H}(t)$ are functions
of time with $N$ components. Mathematically, equation \eqref{eq:yh} can be
replaced by $\bm{y}_{H}(t)=\bm{M}\bm{h}(\omega t)$, where $\bm{M}$
is now a $N \times (2H+1)$ matrix of Fourier coefficients, with each row
corresponding to one of the $N$ components of the $\bm{y}_{H}(t)$. In this case, $\bm{h}(\omega t)$ remains unaltered, and is still given by eqn \eqref{eq:h}. The residual in equation \eqref{eq:f_and_residual} becomes $\bm{r}_{H}(t)$, a function of time with $N$ components, and $\bm{R}$ becomes a $N\times(2H+1)$ matrix. The HB solution is still obtained by solving a nonlinear algebraic system, 
\begin{equation}
\bm{R}(\bm{M})=\bm{0},
\end{equation}
now, with $N(2H+1)$ unknowns.

This describes the general methodology of HB. We now
specialize it for differential constitutive equations in oscillatory
shear, where $\bm{f}(\bm{y},\bm{\dot{y}},t)=\bm{0}$ typically takes the form, 
\begin{equation}
\dot{\bm{y}}+\bm{K}\bm{y}+\bm{f}_{nl}(\bm{y},t)-\bm{f}_{ex}(\bm{y},t)=\bm{0}.\label{eq:general_DE_form}
\end{equation}
Here, $\bm{K}$ is a $N\times N$ diagonal matrix with constant coefficients,
$\bm{f}_{nl}$ is a $N\times 1$ vector of nonlinear functions of the dependent variable $\bm{y}$,
and $\bm{f}_{ex}$ contains the external oscillatory forcing function.
Substituting equation (\ref{eq:yh}) into (\ref{eq:general_DE_form}), and using
equations (\ref{eq:f_and_residual}) and (\ref{eq:residual}), we obtain, 
\begin{equation}
\bm{R}(\bm{M})=\bm{M}\bm{D}+\bm{K}\bm{M}+\bm{F}_{nl}(\bm{M})-\bm{F}_{ex}\left(\bm{M}\right)=\bm{0},\label{eq:residual_eqn}
\end{equation}
where $\bm{F}_{nl}$ and $\bm{F}_{ex}$ are $N \times (2H+1)$ matrices
comprising the Fourier coefficients of the nonlinear $\bm{f}_{nl}$
and external forcing functions $\bm{f}_{ex}$, respectively. $\bm{D}=\text{diag}(-i\omega H,\cdots,0,\cdots,i\omega H)$ is the $(2H+1) \times (2H+1)$ diagonal matrix, alluded to earlier. This  involves solving an algebraic system in $N\times(2H+1)$ unknown elements of $\bm{M}$.

HB is most effective when both the desired solution and the
residual are analytic, and can be expressed in terms of a convergent
Fourier series. When a Fourier series is convergent, Fourier coefficients
decay rapidly with the harmonic index $H$. When the system
is analytic, the $k$th coefficient of the Fourier series decays exponentially
$Y_{k}\sim c^{-k}$, for some constant $c$. In such scenarios, the
difference between the ansatz solution $y_{H}(t)$ and the exact solution
$y(t)$ decays exponentially with $H$ \cite{krack2019harmonic}, 
\begin{equation}
\xi=\left|\left|y_{H}(t)-y(t)\right|\right|_{\infty}\leq u\,e^{-mH},\label{eq:error}
\end{equation}
where $u$ and $m$ are unknown constants, and the infinity norm $\left|\left|x\right|\right|_{\infty}$
denotes the maximum absolute value of $x(t)$ for $t\in\left[0,T\right]$.
Thus, HB shows exponential convergence towards the true solution as
$H$ increases.

\subsection{Harmonic Balance for Giesekus Model}

In this work, we set solvent viscosity $\eta_{s}=0$, so that $\bm{\sigma} = \bm{\sigma}_p$, and the subscript `$p$' can be dropped. For homogeneous flow, the single-mode Giesekus model described by eqns \eqref{eq:totalStress} and \eqref{eq:tensor_form} can be expressed as a set of nonlinear ODEs, 
\begin{equation}
\begin{gathered}\frac{d\sigma_{11}}{dt}+\frac{\sigma_{11}}{\lambda}+\frac{\alpha}{\lambda G}\left(\sigma_{11}^{2}+\sigma_{12}^{2}\right)-2\dot{\gamma}\sigma_{12}=0\\
\frac{d\sigma_{22}}{dt}+\frac{\sigma_{22}}{\lambda}+\frac{\alpha}{\lambda G}\left(\sigma_{22}^{2}+\sigma_{12}^{2}\right)=0\\
\frac{d\sigma_{33}}{dt}+\frac{\sigma_{33}}{\lambda}+\frac{\alpha}{\lambda G}\sigma_{33}^{2}=0\\
\frac{d\sigma_{12}}{dt}+\frac{\sigma_{12}}{\lambda}+\frac{\alpha}{\lambda G}\left(\sigma_{11}+\sigma_{22}\right)\sigma_{12}-\sigma_{22}\dot{\gamma}-G\dot{\gamma}=0.
\end{gathered}
\label{eq:giesekus_eqn}
\end{equation}
Here, $\dot{\gamma}(t)=\gamma_{0}\omega\cos\omega t$ is the oscillatory shear rate, and $\alpha,\lambda$, and $G$ are system parameters. Since the evolution of $\sigma_{33}$ does not depend on other components of the stress tensor $\bm{\sigma}$, we can set the steady state $\sigma_{33}=0$, and retain only $N=3$ dynamic variables, $\bm{y}=[\sigma_{11},\sigma_{22},\sigma_{12}]^{T}$. Thus, equation \eqref{eq:giesekus_eqn} can be cast in the form of equation \eqref{eq:general_DE_form} with a $3\times3$ diagonal matrix $\bm{K}=\lambda^{-1}\bm{I}$, where $\bm{I}$ is the identity matrix. Then $\bm{f}_{nl}$ contains all the quadratic nonlinear terms, 
\begin{equation}
\bm{f}_{nl}=\frac{\alpha}{\lambda G}\left[\begin{array}{c}
(\sigma_{11}^{2}+\sigma_{12}^{2})\\
(\sigma_{22}^{2}+\sigma_{12}^{2})\\
(\sigma_{11}+\sigma_{22})\sigma_{12}
\end{array}\right],\label{eq: fnl}
\end{equation}
and $\bm{f}_{ex}$ consists of forcing terms with $\dot{\gamma}=\gamma_{0}\omega\cos(\omega t)=\frac{\gamma_{0}\omega}{2}\left(e^{i\omega t}+e^{-i\omega t}\right)$,
\begin{equation}
\bm{f}_{ex}=\frac{\gamma_{0}\omega}{2}\left(e^{i\omega t}+e^{-i\omega t}\right)\left[\begin{array}{c}
2\sigma_{12}\\
0\\
(G+\sigma_{22})
\end{array}\right].\label{fex}
\end{equation}
This system can be transformed to the Fourier domain, and solved using HB according to equations \eqref{eq:yh} -- \eqref{eq:residual_eqn}. Using arguments of symmetry, the shear and normal stresses can be represented as truncated Fourier series with only odd and even harmonics, respectively. This reduces the number of unknown coefficients by nearly half.

Consider a Fourier series with $H$ odd harmonics in shear stress. We adopt a convention in which $(2H-1)\omega$ and $(2H-2)\omega$ are the highest harmonics present in shear and normal stresses, respectively. Note that it is possible to choose a convention where $2H$ is the highest harmonic present in the normal stresses. For small values of $\gamma_{0}$ in the linear viscoelastic
regime, this choice is optimal since $H=1$ is sufficient to describe the second harmonic that is excited in $N_1$ and $N_2$.  On the other hand, the MAOS response is typically represented by truncating the power series representation beyond the third harmonic (see eqns \ref{eq:shear_power_series} -- \ref{eq:n2_power_series}), in which case the convention we adopt is optimal since $H=2$ is sufficient. Regardless, the objective of this work is to explore the LAOS regime which calls for modestly large values of $H$. In such scenarios, the specific convention adopted becomes unimportant. Using our convention, the shear and normal stresses are written as, 
\begin{align}
\sigma_{11}(t) & =\sum_{k=-(H-1)}^{H-1}A_{2k}\,e^{i2k\omega t}\nonumber \\
\sigma_{22}(t) & =\sum_{k=-(H-1)}^{H-1}B_{2k}\,e^{i2k\omega t}\nonumber \\
\sigma_{12}(t) & =\sum_{k=-H}^{H-1}C_{2k+1}\,e^{i(2k+1)\omega t}.\label{eq:stress_fourier_series}
\end{align}
There are $2H-1$ Fourier coefficients corresponding to the first
and second normal stresses. Similarly there are $2H$ Fourier coefficients
corresponding to the shear stress. Let 
\begin{equation}
\bm{M}=[\{A_{2k}\}_{k=-(H-1)}^{H-1},\{B_{2k}\}_{k=-(H-1)}^{H-1},\{C_{2k+1}\}_{k=-H}^{H-1}]\label{eq:vectorM}
\end{equation}
 be the vector of unknown coefficients. The total number of
elements in $\bm{M}$ is $(6H-2)$. In principle, these Fourier coefficients
are complex, and the total number of unknowns is twice as many. However,
since stresses are real, Fourier coefficients are complex conjugates;
for example, $A_{2k}=\bar{A}_{-2k}$, where the overbar is used to
denote the conjugate. This means that the real and imaginary parts
are related as $\text{Re}(A_{2k})=\text{Re}(\bar{A}_{-2k})$, and
$\text{Im}(A_{2k})=-\text{Im}(\bar{A}_{-2k})$.

We can substitute eqn \eqref{eq:stress_fourier_series} into
eqn \eqref{eq:giesekus_eqn}, and truncate the quadratic nonlinear term after
the highest harmonic ($2H-1$) is resolved.  For step-by-step derivation of the folowing equations (eqns \eqref{eq:eq_for_T11} -- \eqref{eq:eq_for_T12}) please refer to the supplementary material. For the ODE involving $\sigma_{11}$, we obtain a set of $2H-1$ algebraic nonlinear equations in the Fourier coefficients for $k\in[-(H-1),H-1]$, 
\begin{align}
R_{k}^{11}(\bm{M})= \frac{\alpha}{\lambda G} & \left(\sum_{m=-(H-1)}^{H-1}A_{2m}A_{2(k-m)}+\sum_{m=-H}^{H-1}C_{2m+1}C_{2(k-m)-1}\right) + \nonumber \\
 & \left(i2k\omega+\frac{1}{\lambda}\right)A_{2k} - \gamma_{0}\omega(C_{2k-1}+C_{2k+1})=0.\label{eq:eq_for_T11}
\end{align}
Here and subsequently, all harmonics higher than those prescribed by
the ansatz are assumed to be zero. Thus, $A_{j}=B_{j}=0$ for $\left|j\right|>2H-2$
and $C_{j}=0$ for $\left|j\right|>2H-1$. The quadratic nonlinearity
of the Giesekus model is reflected in the quadratic terms multiplying
the factor $\alpha/\lambda G$. These quadratic terms
are multiplications of two Fourier series and can be determined using
the convolution theorem \cite{krack2019harmonic}. The ODE involving
$\sigma_{22}$ also reduces to a set of $2H-1$ quadratic equations
for $k\in[-(H-1),H-1]$, 
\begin{equation}
\begin{split}
R_{k}^{22}(\bm{M})= \frac{\alpha}{\lambda G} & \left(\sum_{m=-(H-1)}^{H-1}B_{2m}B_{2(k-m)}+\sum_{m=-H}^{H-1}C_{2m+1}C_{2(k-m)-1}\right) + \\
& \left(i2k\omega+\frac{1}{\lambda}\right)B_{2k} = 0.
\end{split}
\label{eq:eq_for_T22}
\end{equation}
Finally, the ODE for $\sigma_{12}$ leads to $2H$ equations
for $k\in[-H,H-1]$ 
\begin{align}
R_{k}^{12}(\bm{M})= & \left(i\left(2k+1\right)\omega+\frac{1}{\lambda}\right)C_{2k+1}  +\frac{\alpha}{\lambda G}\sum_{m=-H}^{H-1}C_{2m+1}\left(A_{2(k-m)}+B_{2(k-m)}\right)\nonumber \\
 & -\gamma_{0}\omega\left(\frac{B_{2k}+B_{2k+2}}{2}\right)-\dfrac{G\gamma_{0}\omega}{2}\delta_{|2k+1|1}=0,\label{eq:eq_for_T12}
\end{align}
where $\delta_{|2k+1|1}$ in the last term of this equation is a Kronecker delta function that is nonzero only when $k=-1$ or $k=0$. This term arises from the externally applied shear force, in which only the first harmonic is triggered.

The residuals corresponding to the equations \eqref{eq:eq_for_T11}
-- \eqref{eq:eq_for_T12} are stacked into an array, 
\begin{equation}
\bm{R}(\bm{M)}=[\{R_{k}^{11}(\bm{M})\}_{k=-(H-1)}^{H-1},\{R_{k}^{22}(\bm{M})\}_{k=-(H-1)}^{H-1},\{R_{k}^{12}(\bm{M})\}_{k=-H}^{H-1}]\label{eq:vectorR}
\end{equation}
 with $6H-2$ elements. We can thus solve these nonlinear equations
for $\bm{M}$ by using a suitable nonlinear solver that seeks to minimize
the residuals. Once we obtain the solution $\bm{M}$, and therefore all the
Fourier coefficients, we can compute the dynamic moduli used in the
Fourier representation, eqns \eqref{eq:shear} -- \eqref{eq:secondNormal},
via, $G_{n}' =i\left(C_{n}-C_{-n}\right)$, $G_{n}''  =\left(C_{n}+C_{-n}\right)$,
$F_{n}'  =i\left(\left(A_{n}-B_{n}\right)-\left(A_{-n}-B_{-n}\right)\right)$,  $F_{n}'' =\left(\left(A_{n}-B_{n}\right)+\left(A_{-n}-B_{-n}\right)\right)$,
$S_{n}'  = i\left(B_{n}-B_{-n}\right)$, and $S_{n}'' =\left(B_{n}+B_{-n}\right)$.

\subsection{Solving Harmonic Balance Equations}

HB leads to a system of nonlinear equations in $6H-2$ unknowns
stacked in the vector $\bm{M}$ (eqn. \eqref{eq:vectorM}), which
can be expressed in terms of the $6H-2$ elements of the residual
vector $\bm{R}$ (eqn. \eqref{eq:vectorR}). We solve these equations
using the MATLAB function \texttt{fsolve}. It uses
a trust region dogleg algorithm, which belongs to a popular class
of methods for optimization problems. It attempts to solve for the
elements of $\bm{M}$ by minimizing the elements of the residual $\bm{R}$.
We set the termination criterion as $||\bm{R}||_{\infty}<10^{-12}$.

In this work, the LAOS response of the Giesekus model is
computed at different frequencies $\omega$, and three decades of strain
amplitude $\gamma_{0}$ spanning the range $10^{-2}-10^{1}$. The
input to \texttt{fsolve} is the residual vector $\bm{R}(\bm{M})$ as a function of the unknown coefficients $\bm{M}$, and an initial estimate of the solution vector. For small strains, $\gamma_{0}\lesssim0.1$, we use the analytical MAOS response as the initial guess \cite{kate2012large}. All Fourier coefficients corresponding to harmonics greater than the third harmonic are set to zero. This
is a reasonably good initial guess for small strains. This is important
because nonlinear systems can have multiple solutions, some of which
may even be unstable. A good initial guess mitigates the risk of converging
to unstable solutions, as \texttt{fsolve} typically
converges to the solution closest to it.

As $\gamma_{0}$ increases beyond 0.1, the solution gradually
begins to deviate from the MAOS response. Hence, we use a laddered
approach to develop good initial guesses for intermediate and large
values of $\gamma_{0}$. In this work, the ladder comprises three
intermediate rungs at strain amplitudes given by $[1.0,2.2,4.6]$
that are dispersed between 0.1 and 10. Thus, we selected $\gamma_{0}=0.1$
and $\gamma_{0}=1$ as the first two steps of this ladder, and two intermediate rungs between $\gamma_{0}=1-10$, that are logarithmically equispaced. For any given value of $\gamma_{0}$ between $0.1-10$, say $\gamma_{0}=3.5$, we first use the analytical
MAOS solution as an initial estimate for the response at $\gamma_{0}=0.1$.
We then use the HB solution obtained at $\gamma_{0}=0.1$ as the initial
guess for the $\gamma_{0}$ corresponding to the next rung in the
ladder, which is $\gamma_{0}=1$. This process is repeated until the
value of $\gamma_{0}$ at which the solution is desired is smaller
than the $\gamma_{0}$ that corresponds to the next rung of the ladder.
By staging the solution in this manner, we ensure that a good initial
guess is available for the solver. This gradually activates higher
harmonics as $\gamma_{0}$ is increased, and avoids convergence to
incorrect solution branches, especially at high frequencies.

\subsection{Solving IVP}

To solve the IVP, we first recast the ODEs in eqn \eqref{eq:giesekus_eqn}
in dimensionless form. If we do not non-dimensionalize the ODEs, the numerical
solution becomes unstable at large $\omega$ and $\gamma_{0}$. Using
$\tilde{t}=t/\lambda$, $\tilde{\bm{\sigma}}=\bm{\sigma}/(G\lambda\omega\gamma_{0})$,
$\tilde{\dot{\gamma}}=\dot{\gamma}/(\gamma_{0}\omega)=\cos(\De\,\tilde{t})$,
$\De=\lambda\omega$, and Wiessenberg number $\Wi=\lambda\gamma_{0}\omega$ \cite{khair2016large}, we obtain, 
\begin{align}
r_{11}(\tilde{t}) & =\frac{d\tilde{\sigma}_{11}}{d\tilde{t}}+\tilde{\sigma}_{11}+\alpha\Wi\left(\tilde{\sigma}_{11}^{2}+\tilde{\sigma}_{12}^{2}\right)-2\tilde{\dot{\gamma}}\tilde{\sigma}_{12}\Wi \approx 0\nonumber \\
r_{22}(\tilde{t}) & =\frac{d\tilde{\sigma}_{22}}{d\tilde{t}}+\tilde{\sigma}_{22}+\alpha\Wi\left(\tilde{\sigma}_{22}^{2}+\tilde{\sigma}_{12}^{2}\right) \approx 0\label{eq:dim_giesekus}\\
r_{12}(\tilde{t}) & =\frac{d\tilde{\sigma}_{12}}{d\tilde{t}}+\tilde{\sigma}_{12}+\alpha\Wi\left(\tilde{\sigma}_{11}+\tilde{\sigma}_{22}\right)\tilde{\sigma}_{12}-\tilde{\dot{\gamma}}-\tilde{\dot{\gamma}}\tilde{\sigma}_{22}\Wi \approx 0\nonumber 
\end{align}
The residuals $r_{11}(\tilde{t}),r_{22}(\tilde{t})$, and $r_{12}(\tilde{t})$
should ideally be zero if the exact solution is inserted in the above
equations. However, numerical approximation results in nonzero residuals,
which are nevertheless themselves periodic.

We solve the ODEs given by eqn \eqref{eq:dim_giesekus} using
the adaptive implicit Runge-Kutta method \texttt{ode23s}
in MATLAB. In order to obtain the alternance solution, we monitor
the relative difference between successive peaks of the stress outputs.
This difference falls as the initial transients decay, and eventually
becomes equal to zero. Once the relative difference falls below a
suitable threshold ($10^{-6}$ in this work), we assume that the periodic
steady state has been attained. The number of oscillation cycles or
time required to reach steady state depends on the magnitude of $\gamma_{0}$
and $\omega$. Once the periodic steady state solution is obtained,
it is interpolated at $N_{t}=1000$ uniformly spaced grid points
per cycle using cubic splines. Fourier representations of the numerically computed stresses are obtained using fast Fourier transform (FFT). This allows us to write,
\begin{equation}
\bm{\tilde{\sigma}}_{\text{IVP}}(\tilde{t})=\left[\begin{array}{c}
\tilde{\sigma}_{11}\\
\tilde{\sigma}_{22}\\
\tilde{\sigma}_{12}
\end{array}\right]=\frac{1}{G\Wi}\left[\begin{array}{c}
\sum_{n}\gamma_{0}^{2}\left(P_{n}'\sin n\,\De\,\tilde{t}+P_{n}''\cos n\,\De\,\tilde{t}\right)\\
\sum_{n}\gamma_{0}^{2}\left(Q_{n}'\sin n\,\De\,\tilde{t}+Q_{n}''\cos n\,\De\,\tilde{t}\right)\\
\sum_{n}\gamma_{0}\left(G_{n}'\sin n\,\De\,\tilde{t}+G_{n}''\cos n\,\De\,\tilde{t}\right)
\end{array}\right],\label{eq:sigma}
\end{equation}
where the summation over $n$ extends to the largest nonzero mode,
and includes only even and odd harmonics for normal and shear stresses,
respectively. The normal stress coefficients $P_{n}'$ and $P_{n}''$
(and $Q_{n}'$ and $Q_{n}''$) are analogous to $G_{n}'$ and $G_{n}''$
for the shear stress. They can be used to obtain the coefficients
corresponding to the first and second normal stress differences, 
\begin{align}
F_{n}' & =P_{n}'-Q_{n}', & F_{n}'' & =P_{n}''-Q_{n}'',\nonumber \\
S_{n}' & =Q_{n}', & S_{n}'' & =Q_{n}''.
\end{align}

\subsection{Comparing HB and NI solutions}

How can we compare the accuracy of the solution obtained
using HB and NI, especially in the LAOS regime, where the exact solution is not known? One approach, and the one we adopt here, is to compare the magnitude of the residual or error obtained by inserting the numerical solution back into the set of ODEs. We describe how to compute this error term, first for NI, and then for HB.

Once the numerical solution $\bm{\tilde{\sigma}}_{\text{IVP}}(\tilde{t})$
is obtained in the form shown in equation \eqref{eq:sigma}, we can
symbolically compute the derivative at any arbitrary time as, 
\begin{equation}
\frac{d\bm{\tilde{\sigma}}_{\text{IVP}}}{d\tilde{t}}=\frac{1}{G\Wi}\left[\begin{array}{c}
\sum_{n}n\De\gamma_{0}^{2}\left(P_{n}'\cos\left(n\,\De\,\tilde{t}\right)-P_{n}''\sin\left(n\,\De\,\tilde{t}\right)\right)\\
\sum_{n}n\De\gamma_{0}^{2}\left(Q_{n}'\cos\left(n\,\De\,\tilde{t}\right)-Q_{n}''\sin\left(n\,\De\,\tilde{t}\right)\right)\\
\sum_{n}n\De\gamma_{0}\left(G_{n}'\cos\left(n\,\De\,\tilde{t}\right)-G_{n}''\sin\left(n\,\De\,\tilde{t}\right)\right)
\end{array}\right],\label{eq:dSigma}
\end{equation}
where the summation over $n$ runs over indices identical to equation
\eqref{eq:sigma}. We can substitute equations \eqref{eq:sigma} and
\eqref{eq:dSigma} in the dimensionless ODEs (eqn \eqref{eq:dim_giesekus}),
and compute the residuals $r_{11}(\tilde{t})$, $r_{22}(\tilde{t})$,
and $r_{12}(\tilde{t})$.

To compare solutions over a period of oscillation, we divide
a single dimensionless period $\tilde{T}=2\pi/\De$ into $N_{t}=1000$
uniformly spaced points $\tilde{t}_{i}=i\Delta\tilde{t}$ for $i=0,\cdots,N_{t}-1$,
with $\Delta\tilde{t}=\tilde{T}/N_{t}$. Residuals are evaluated at
these points, and stacked into a vector, $\bm{r}=[r_{11}(\tilde{t}_{i}),r_{22}(\tilde{t}_{i}),r_{12}(\tilde{t}_{i})]_{i=0}^{N_{t}-1}$
with $3N_{t}$ elements. The error associated with this term can be
computed from the root mean square value of $\bm{r}$ via the Euclidean
or 2-norm as, 
\begin{equation}
\epsilon_{r}=\dfrac{1}{3N_{t}}||\bm{r}||_{2}.\label{eq:residual2}
\end{equation}

We can follow a nearly identical process to compute $\epsilon_{r}$
from the HB solution $\bm{\sigma}_\text{HB}$ which is obtained in a form
given by eqn \eqref{eq:stress_fourier_series}. We then (i) non-dimensionalize
this solution similar to eqn \eqref{eq:sigma}, (ii) symbolically
compute the first derivative of eqn \eqref{eq:stress_fourier_series}
using \eqref{eq:dSigma}, (iii) insert these expressions in the dimensionless
ODEs eqn \eqref{eq:dim_giesekus}, (iv) compute the residual vector
$\bm{r}$ at the same points $\{\tilde{t}_{i}\}$, and (v) evaluate
the error $\epsilon_{r}$.

\section{Results and Discussion}
\label{sec:results}

We examine results for the oscillatory shear response of the Giesekus model obtained using HB. In all the calculations performed in this study, we assume $G=1$ Pa, $\lambda=1$ s, $\alpha=0.3$ and $\eta_{s}=0$ Pa$\cdot$s. We analyze the periodic stress waveforms and harmonic coefficients over a range of $\gamma_0$ and $\omega$. We study the convergence rate as the number of harmonics in the ansatz is increased. Finally, we compare HB and NI in terms of accuracy and computational cost.

\subsection{Stress waveforms and Lissajous plots}

\begin{figure}
\begin{center}
\begin{tabular}{cc}
\includegraphics[scale=0.3]{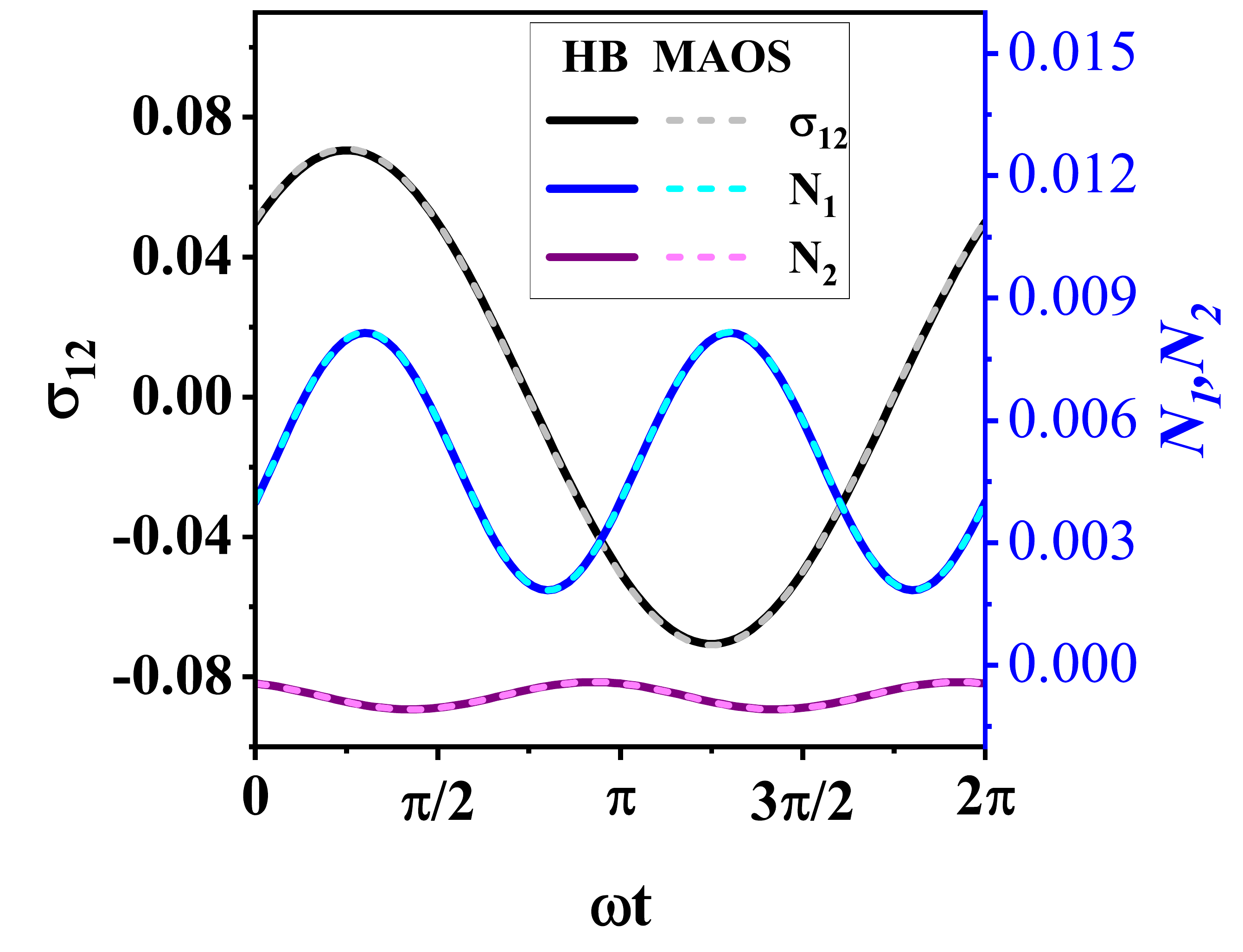}  & \includegraphics[scale=0.085]{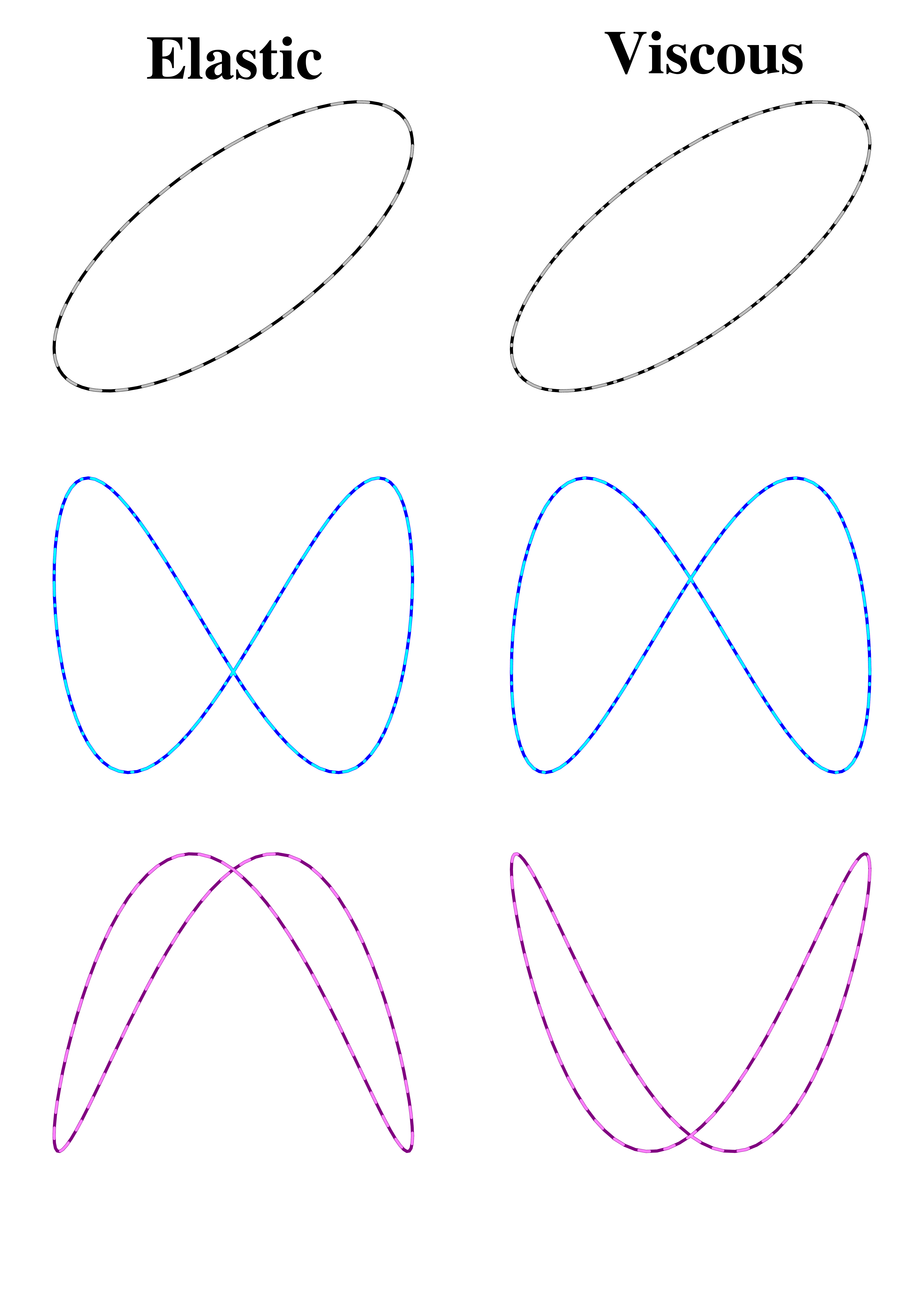}  \\
(a) stress waveform & (b) Lissajous curves \\
\end{tabular}\caption{\label{fig:periodic_data_1} The stress response of the Giesekus model at $\gamma_{0}=0.1$ and $\De=1$ computed with HB method (solid lines) overlaps the MAOS analytical solution (dashed lines) for the (a) shear and normal stress waveforms, and (b) the corresponding elastic and viscous Lissajous curves.}
\end{center}
\end{figure}

We begin by comparing the computed oscillatory stress response using HB with known analytical solutions in the weakly nonlinear or MAOS regime with $\gamma_{0}=0.1$, and $\omega=1$ rad/s. In the HB method, we set $H=5$. The resulting periodic stress waveforms for $\sigma_{12}$,
$N_{1}$, and $N_{2}$ are shown in figure \ref{fig:periodic_data_1}(a).
The waveforms are approximately sinusoidal, and do not exhibit the
effect of strong nonlinearities. The corresponding normalized elastic
and viscous Lissajous curves are shown in figure \ref{fig:periodic_data_1}(b).
The elastic Lissajous curves plot the stress versus strain normalized
by their maximum values. For the shear stress, this shows $\sigma_{12}(t)/\sigma_{12}^{\max}$ versus $\gamma(t)/\gamma_{0}=\sin\omega t$, where $\sigma_{12}^{\max}$ corresponds to the magnitude of the maximum value of $\sigma_{12}(t)$ in figure \ref{fig:periodic_data_1}(a). Similarly, the viscous Lissajous curves plot the stress versus strain rate, normalized by their maximum
values. For $\gamma_{0}=0.1$, the peak magnitude of $\sigma_{12}$ is greater than the peak magnitude of $N_{1}$ and $N_{2}$. $N_{2}$ is smaller than $N_{1}$, and is negative. The analytical MAOS solution for the shear and normal stress differences is overlaid on both plots using dashed lines, and essentially overlaps with the computed solution. 

\begin{figure}
\centering
\begin{tabular}{cc}
\includegraphics[scale=0.3]{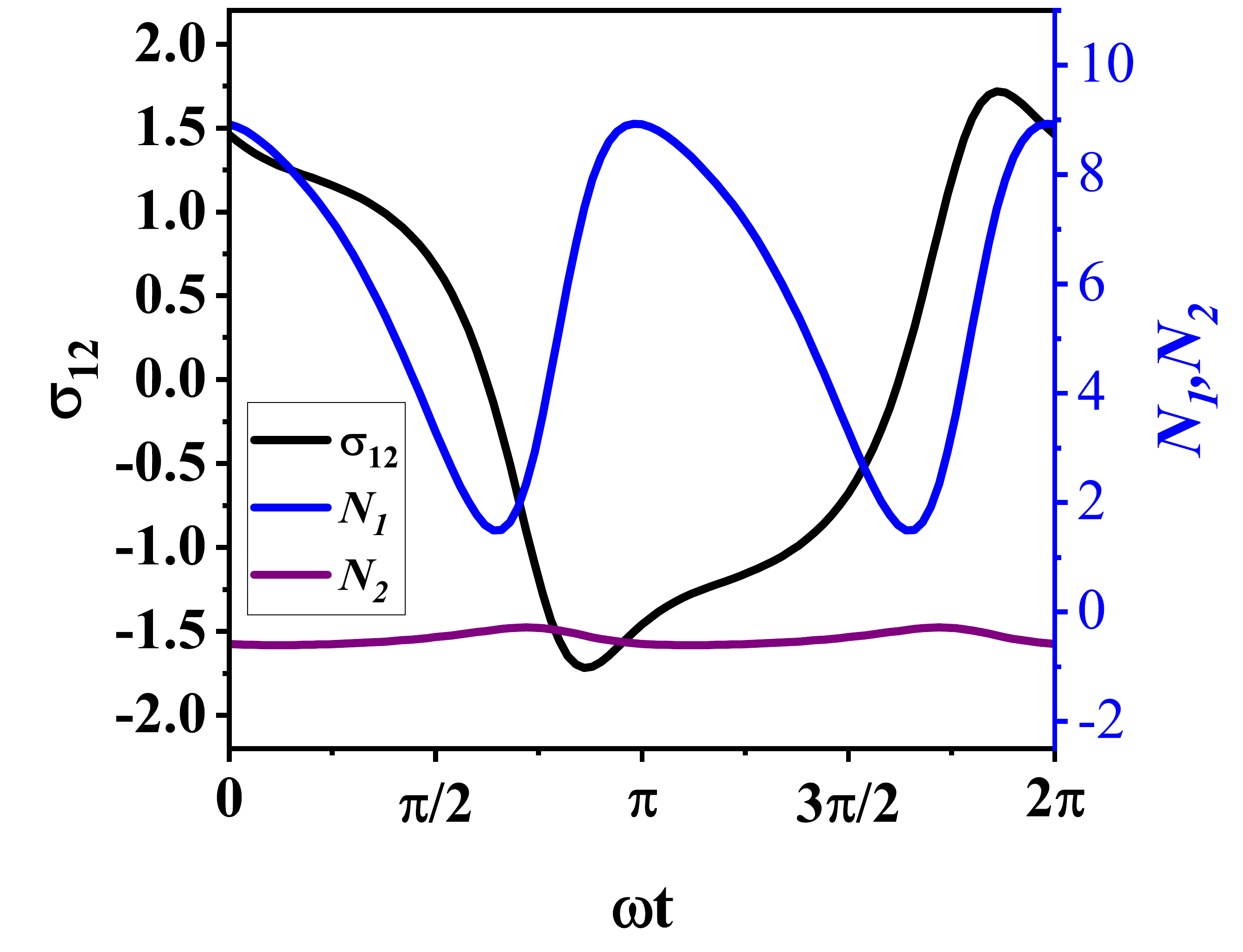}  & \includegraphics[scale=0.082]{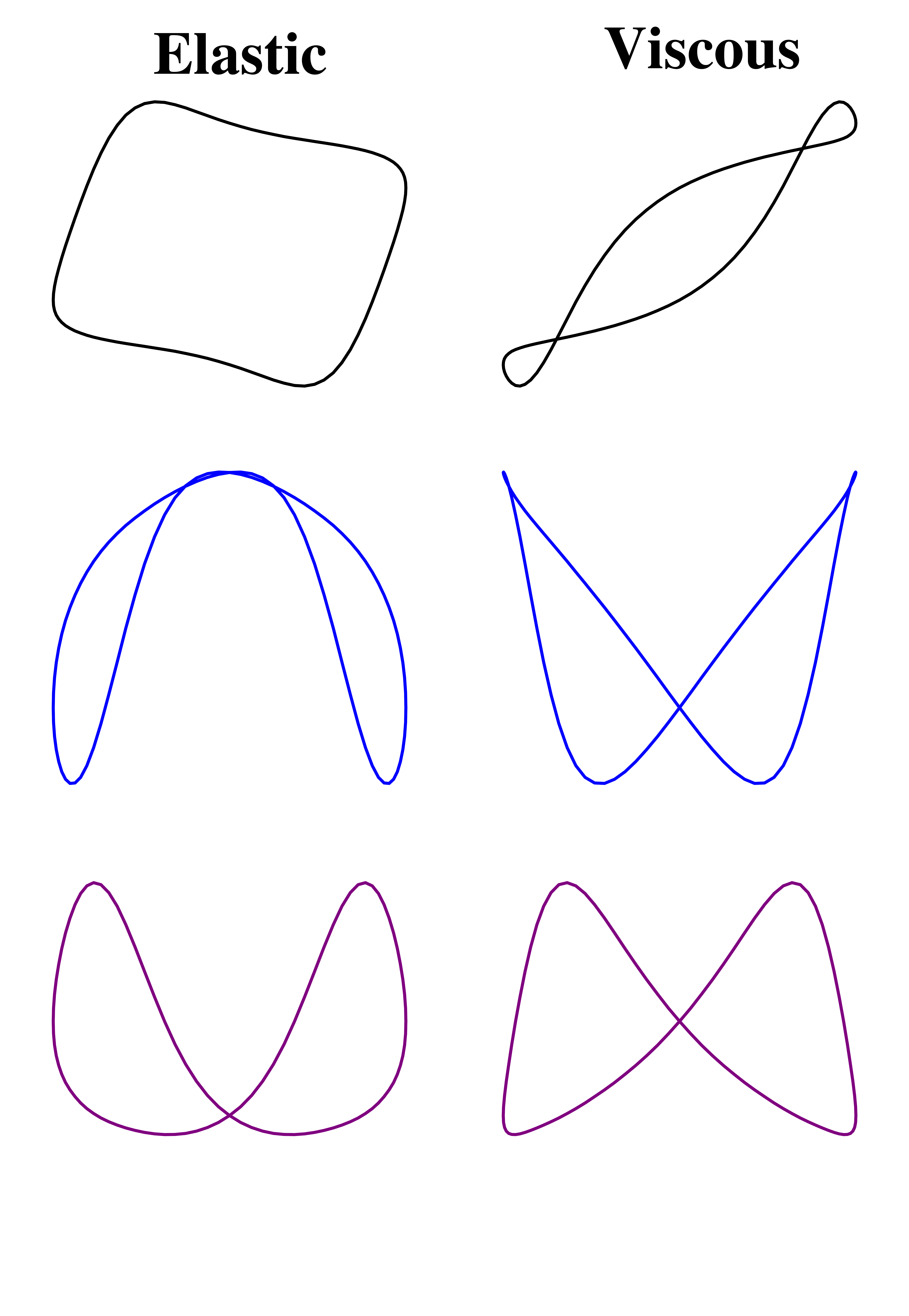}  \\
(a) stress waveform  & (b) Lissajous curves \\
\end{tabular}\caption{\label{fig:periodic_data_2} The oscillatory stress response of the Giesekus model at $\gamma_{0}=10$ and $\De=1$ computed with HB method.
The computed (a) shear and normal stress waveforms, and the corresponding
(b) elastic and viscous Lissajous curves are shown, similar to figure
\ref{fig:periodic_data_1}.}
\end{figure}

Figure \ref{fig:periodic_data_2} shows results similar to
figure \ref{fig:periodic_data_1} at a frequency of $\omega=$ 1 rad/s,
but at a higher strain amplitude of $\gamma_{0}=10$. This corresponds
to $\De=1$ and $\Wi=10$, and
strongly activates the nonlinear terms in the model. The stress waveforms
are significantly distorted, and far from sinusoidal. The peak magnitude
of $N_{1}$ becomes significantly larger than the peak magnitudes
of $\sigma_{12}$ and $N_{2}$. Mathematically, the combination of moderate $\De$ and large $\Wi$ leads to a tightly coupled system which is manifested in the emergence of secondary loops in the Lissajous curves. The occurrence of self-intersections in Lissajous curves has been reported for many experimental systems \cite{jeyaseelan2008network, Ewoldt2008, tee1975nonlinear, ewoldt2010large, khandavalli2016comparison, moud2021viscoelastic, yazar2017non, yasin2021large}, and constitutive models \cite{jeyaseelan2008network, leygue2006tube, hyun2013numerical}, including the Giesekus model \cite{ewoldt2010secondary, kate2012large}.

\begin{figure}
\begin{center}
\begin{tabular}{cc}
\includegraphics[scale=0.3]{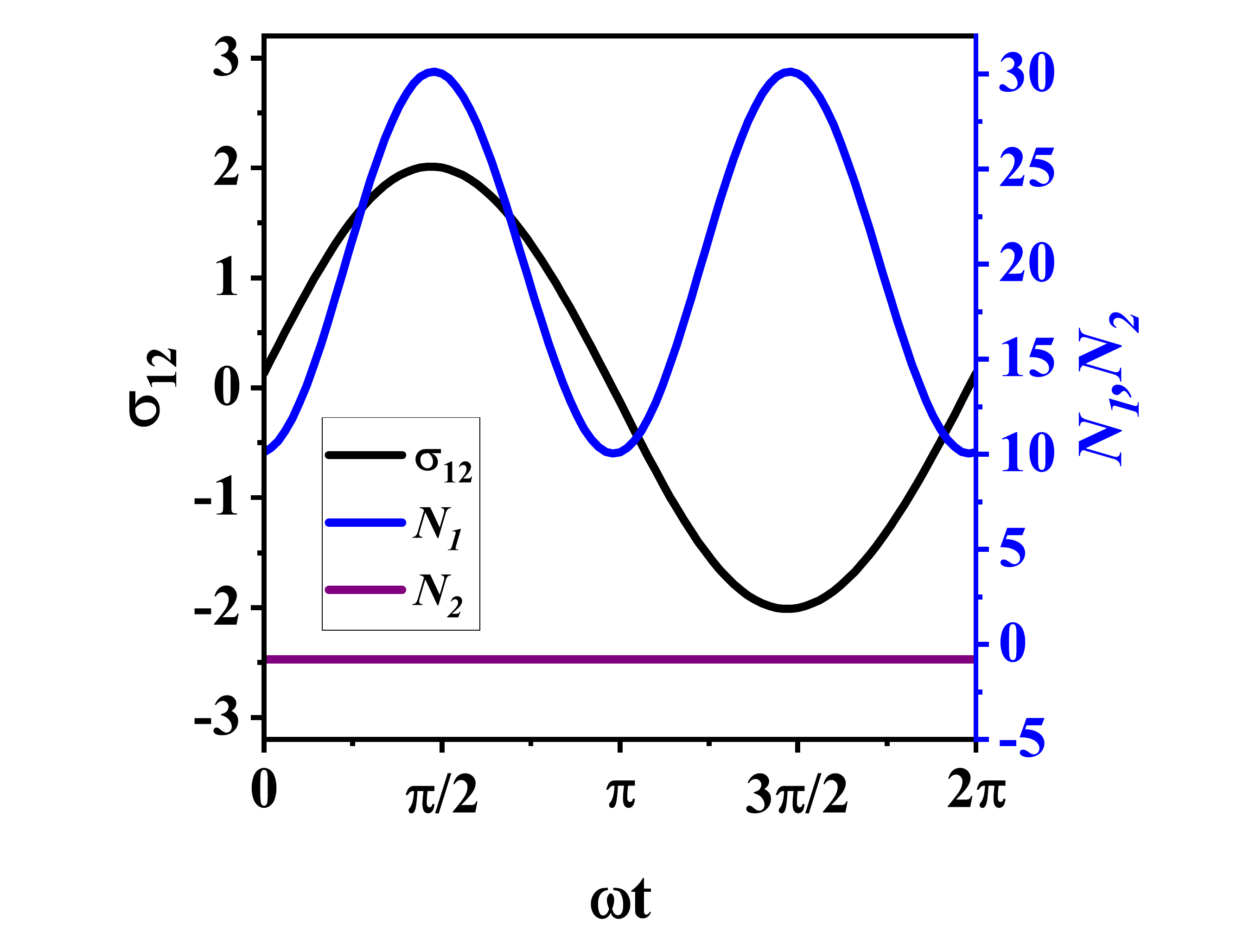} & \includegraphics[scale=0.082]{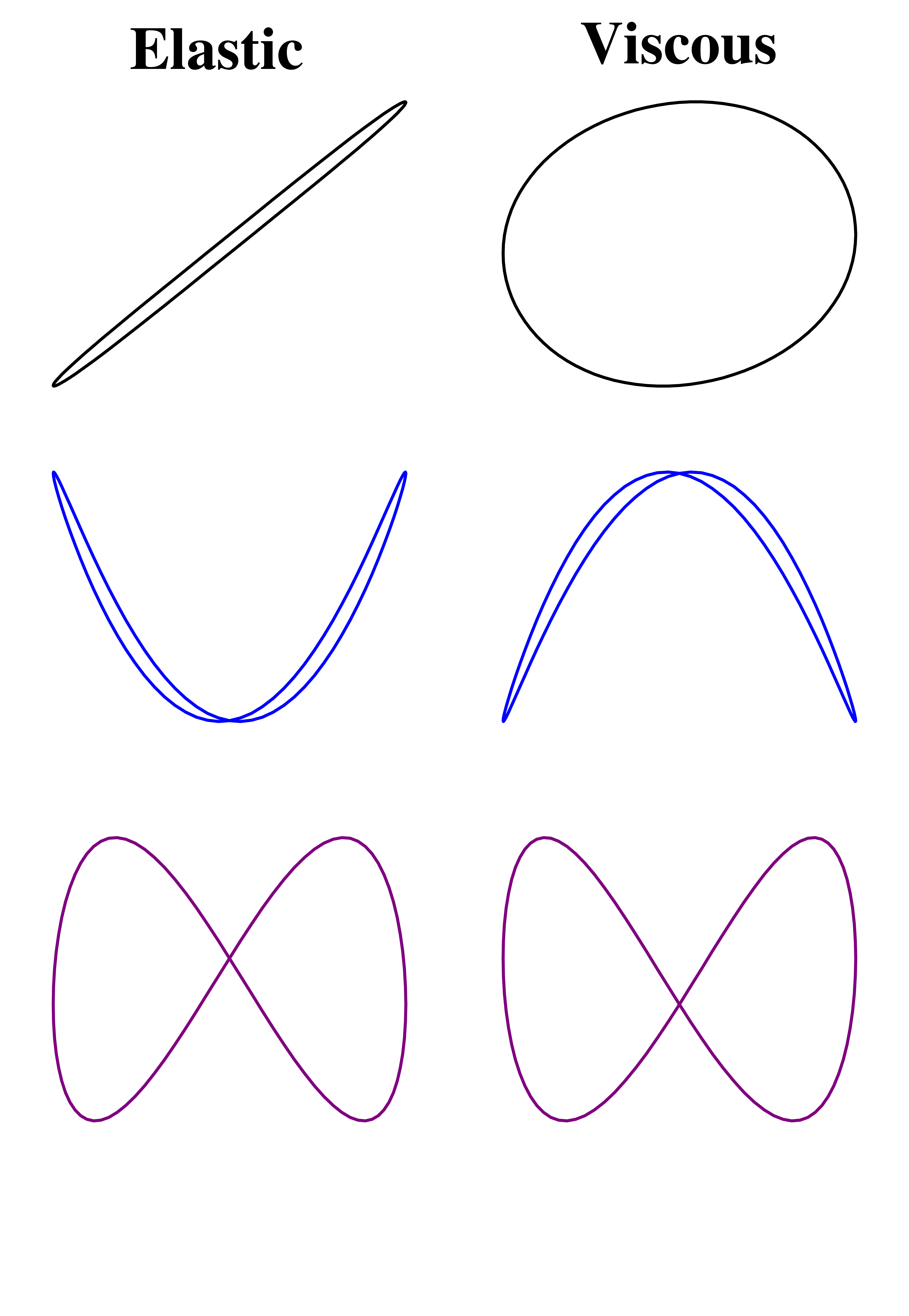} \\
(a) stress waveform & (b) Lissajous curves \\
\end{tabular}
\end{center}
\caption{\label{fig:periodic_data_3} The oscillatory stress response of the Giesekus model at $\gamma_{0} = 10$ and $\De=100$ computed with HB method. The computed (a) shear and normal stress waveforms, and the corresponding (b) elastic and viscous Lissajous curves are shown, similar to figure \ref{fig:periodic_data_1}.}
\end{figure}



Secondary loops disappear at large frequency $\omega=100$ rad/s as shown in figure \ref{fig:periodic_data_3}. Here $\gamma_{0} = 10$, similar to figure \ref{fig:periodic_data_2}. Despite the large $\De=100$ and $\Wi=1000$, the stress waveform shown in \ref{fig:periodic_data_3}(a) is not distorted, and the Lissajous curves no longer exhibit self-intersections associated with a strong viscoelastic response. Elasticity dominates material response at such high frequencies resulting in surprisingly simple profiles.

\begin{figure}
\centering
\includegraphics[scale=0.2]{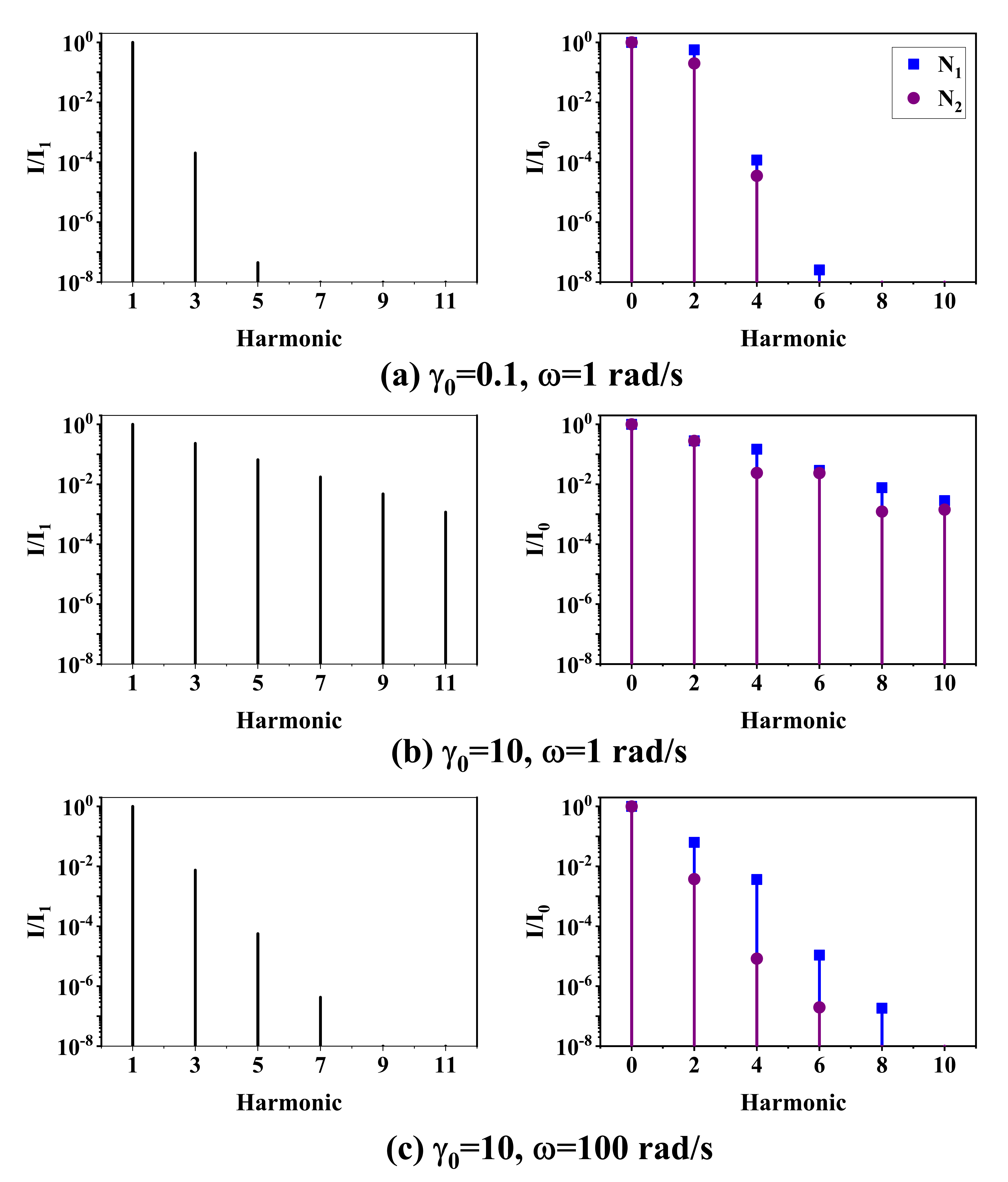} \caption{\label{fig:intensity}Normalized harmonic intensities for shear stress $I/I_{1}$ (left column), and normal stress difference $I/I_{0}$ (right column) at (a) $\gamma_{0}=0.1$, $\omega=1$ rad/s, (b) $\gamma_{0}=10$, $\omega=1$ rad/s, and (c) $\gamma_{0}=10$, $\omega=100$ rad/s using HB.}
\end{figure}

Figure \ref{fig:intensity} plots the intensities of the
harmonics for shear and normal stress differences corresponding to
the three cases studied in figures \ref{fig:periodic_data_1} --
\ref{fig:periodic_data_3}. For shear stress, the intensity of the
$k$th harmonic is defined as $I_{k}=\sqrt{(G_{k}^{\prime})^{2}+(G_{k}^{\prime\prime})^{2}}$.
Analogous expressions are used for normal stress differences. In figure
\ref{fig:intensity}, these intensities are normalized by the leading
harmonic i.e., $I_{1}$ for shear stress, and $I_{0}$ for normal
stress differences. These plots quantitatively reinforce our qualitative
observations from the stress waveforms and Lissajous curves. In figure
\ref{fig:intensity}(a), the weakness of harmonics beyond the first for the shear stress, and beyond the second for normal stresses, suggests that linear viscoelastic modes still dominate the response. On the other hand, higher harmonics exhibit slow decay in figure \ref{fig:intensity}(b) due to the tight coupling in Giesekus model equations, with $I_{3}/I_{1}$ being $23\%$ and $I_{4}/I_{0}$ being $14\%$. The contribution of higher harmonics is small in figure \ref{fig:intensity}(c), despite being deep in the LAOS regime, due to the dominance of elastic contributions. Thus, a strong deformation does not guarantee the activation of larger modes and hence, a judicious choice of $H$ is necessary to ensure convergence of the HB method while minimizing the computational effort.


\begin{figure*}
\centering{}\includegraphics[scale=0.25]{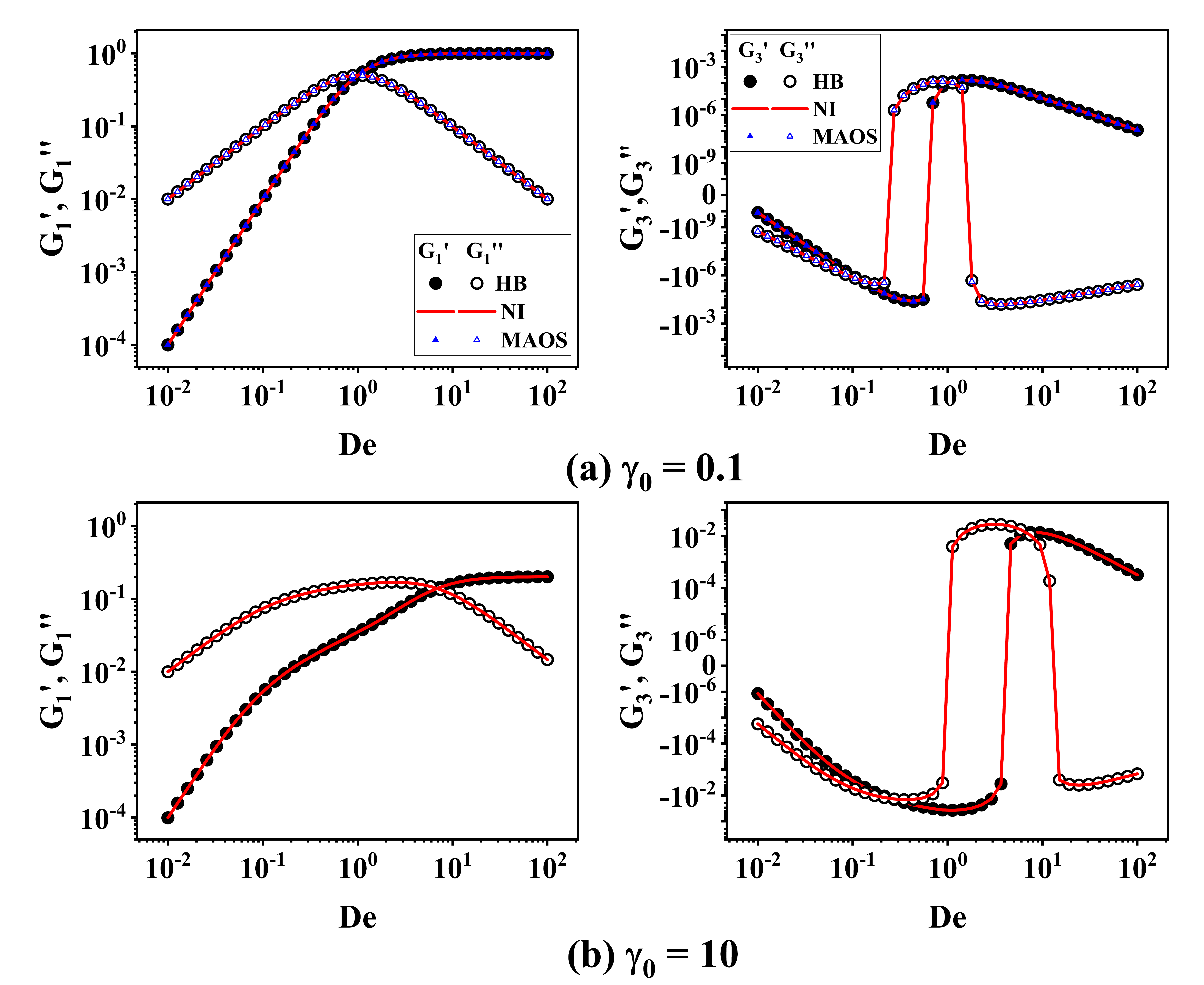} \caption{\label{fig:shear_coeffs} The Fourier coefficients of shear stress
$(G_{1}',G_{1}'',G_{3}',$ and $G_{3}'')$ at (a) $\gamma_{0}=0.1$
and (b) $\gamma_{0}=10$ using HB, and NI. The analytical MAOS solution shown in subfigure
(a) overlaps with the HB solution.}
\end{figure*}

\subsection{Frequency Sweeps of Leading Nonlinear Harmonics}

Figure \ref{fig:shear_coeffs} shows the first and third harmonic coefficients for shear stress at $\gamma_{0}=0.1$ and $\gamma_{0}=10$ across four decades of frequency, $\omega\in[10^{-2},10^{2}]$. For the HB solution, storage moduli are shown by filled circles, while loss moduli are shown by open circles. The analytical MAOS solution is shown by small triangles \cite{kate2012large}. The third harmonic takes both positive and negative values. Hence, we use a symmetrical logscale (or symlog) to represent the vertical axis. In these plots, a small linear region is introduced near the origin to account for sign changes. The HB solution is in good agreement with the analytical
MAOS solutions. 

The nonlinear regime in oscillatory shear can be crudely
defined by $\De>1$ and $\Wi>1$. Figure \ref{fig:shear_coeffs} shows the expected increase in relative intensity of the third harmonic and decrease in the magnitude of primary harmonics as we move towards the LAOS
regime. The third harmonic becomes more than 10\% of the first harmonic at $\gamma_{0}=10$ which leads to a distorted stress waveform and
secondary loops in viscous Lissajous-Bowditch curves in the range of $0.1<\De<3$, as discussed previously. 

\begin{figure}
\begin{center}
\includegraphics[scale=0.25]{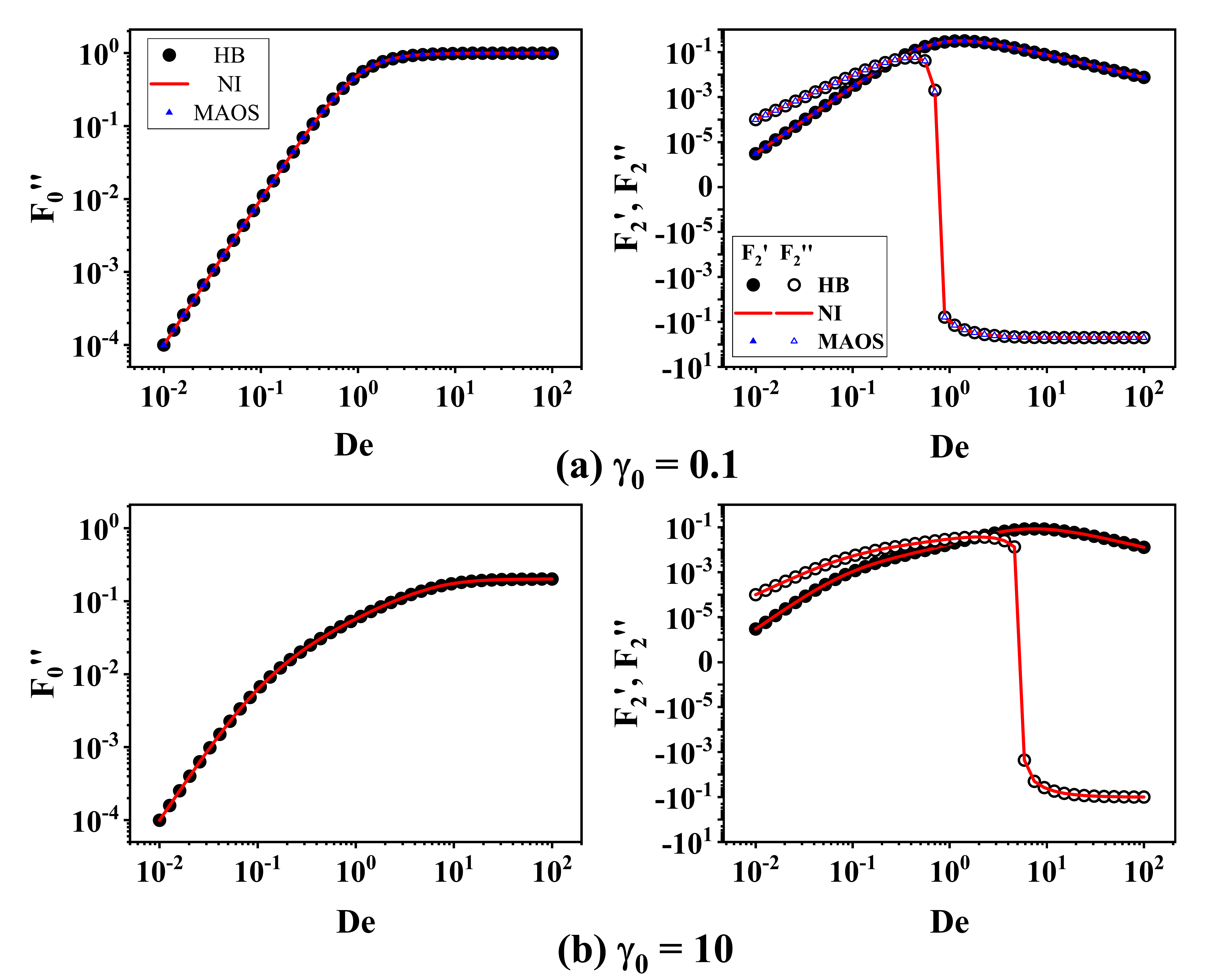}
\par
\end{center}
\caption{\label{fig:n1_coeffs} The Fourier coefficients of the first normal
stress difference $(F_{0}'',F_{2}'\text{ and }F_{2}'')$ at (a) $\gamma_{0}=0.1$,
and (b) $\gamma_{0}=10$ using harmonic balance (HB), and numerical
integration (NI). The analytical MAOS solution is shown in subfigure
(a).}
\end{figure}

\begin{figure}
\begin{center}
\includegraphics[scale=0.25]{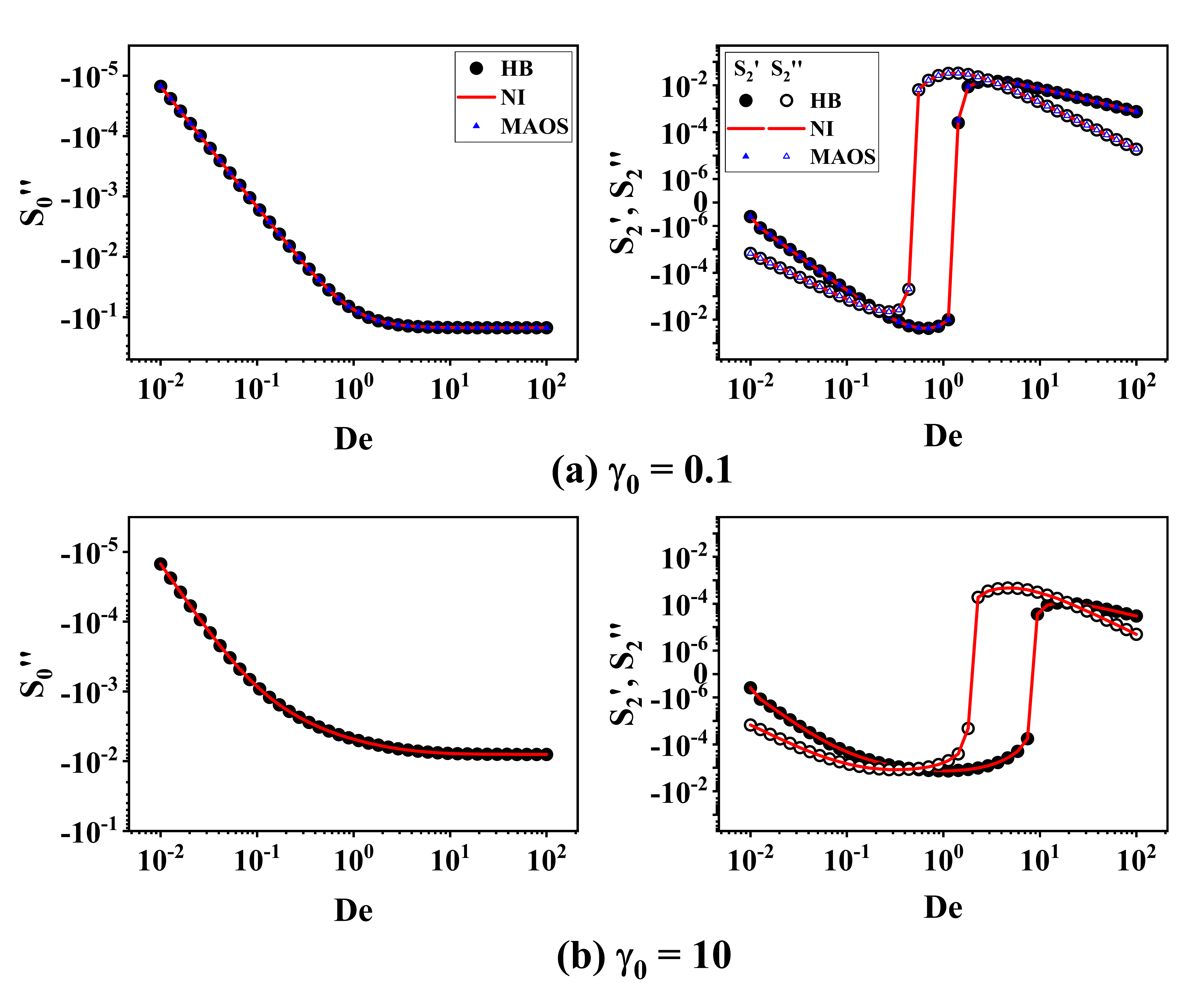}
\par
\end{center}
\caption{\label{fig:n2_coeffs} The Fourier coefficients of the second normal
stress difference $(S_{0}'',S_{2}',$ and $S_{2}'')$ at (a) $\gamma_{0}=0.1$,
and (b) $\gamma_{0}=0.01$ using harmonic balance (HB), and numerical
integration (NI). The analytical MAOS solution is shown in subfigure
(a). }
\end{figure}

Figure \ref{fig:n1_coeffs} and \ref{fig:n2_coeffs} represent the corresponding coefficients for the first normal stress difference and second normal stress difference, respectively, over four decades of frequency at $\gamma_{0}=0.1$ and $\gamma_{0}=10$. The normal stress coefficients are found by normalizing the harmonics with $\gamma_{0}^{2}$. Notice that the leading order terms for $N_{1}$
and $N_{2}$ are always positive and negative respectively. Similar to our observations for the shear stress coefficients, the magnitude of primary harmonics ($0\omega$ and $2\omega$ terms) decline as the contributions of higher frequencies become significant at larger $\gamma_{0}$. In addition, the frequencies at which $F_{2}'',S_{2}'$ and $S_{2}''$ cross the x-axis shift to higher values with increasing $\gamma_{0}$. 

The results at $\gamma_{0}=0.1$ are in accordance with the analytical MAOS solution \cite{kate2012large,nam2008prediction}. In the linear regime, as the system moves from low-frequency viscous response to a high frequency elastic response, $F_{0}''$ and $F_{2}''$ plateau, which is expected since they are functions of $G'$, while $F_{2}'$ peaks and then decreases due to its dependence
on $G''$. This can be understood in terms of MAOS relationships \cite{kate2012large, nam2008prediction}, $F_{0}''=G'(\omega)$, $F_{2}'= G''(\omega) - G''(2\omega)/2$, and $F_{2}'' = -G'(\omega) + G'(2\omega)/2$.

Interestingly, these relations also describe $N_1$ for the UCM model (eqn \eqref{eq:N1_ucm}), which is a special case of the Giesekus model with $\alpha=0$ and $\eta_{s}=0$. Thus, we may say that in the MAOS limit, where only $\omega,2\omega$ and $3\omega$ terms are activated, $G',G'',F_{0}'',F_{2}'$ and $F_{2}''$ are not functions of $\alpha$ and behave essentially like a UCM fluid. However, the distinction between the Giesekus and UCM models becomes apparent from $N_2$, which is zero for the UCM. For the Giesekus model, the Fourier coefficients of $N_{2}$ in the MAOS regime are non-zero and show a dependence on $\alpha$. In addition, these coefficients
$S_{0}'',S_{2}'$ and $S_{2}''$ cannot be expressed as a linear combination of $G'$ and $G''$, as opposed to $F_{0}'',F_{2}'$ and $F_{2}''$. Thus, even in the linear regime, $N_2$ is a function of $\alpha$ and $\omega$, while $\sigma_{12}$ and $N_1$ are dependent only on $\omega$. In principle, we can map a real system to a Giesekus model in the linear limit using $N_{2}$ to determine $\alpha$. However, it is practically difficult since the signal is weak in this regime  \cite{nam2008prediction}. Thus, LAOS analysis of nonlinear stress waveforms is necessary to fit experimental data to the Giesekus model. 

\begin{figure}
\begin{centering}
\includegraphics[scale=0.25]{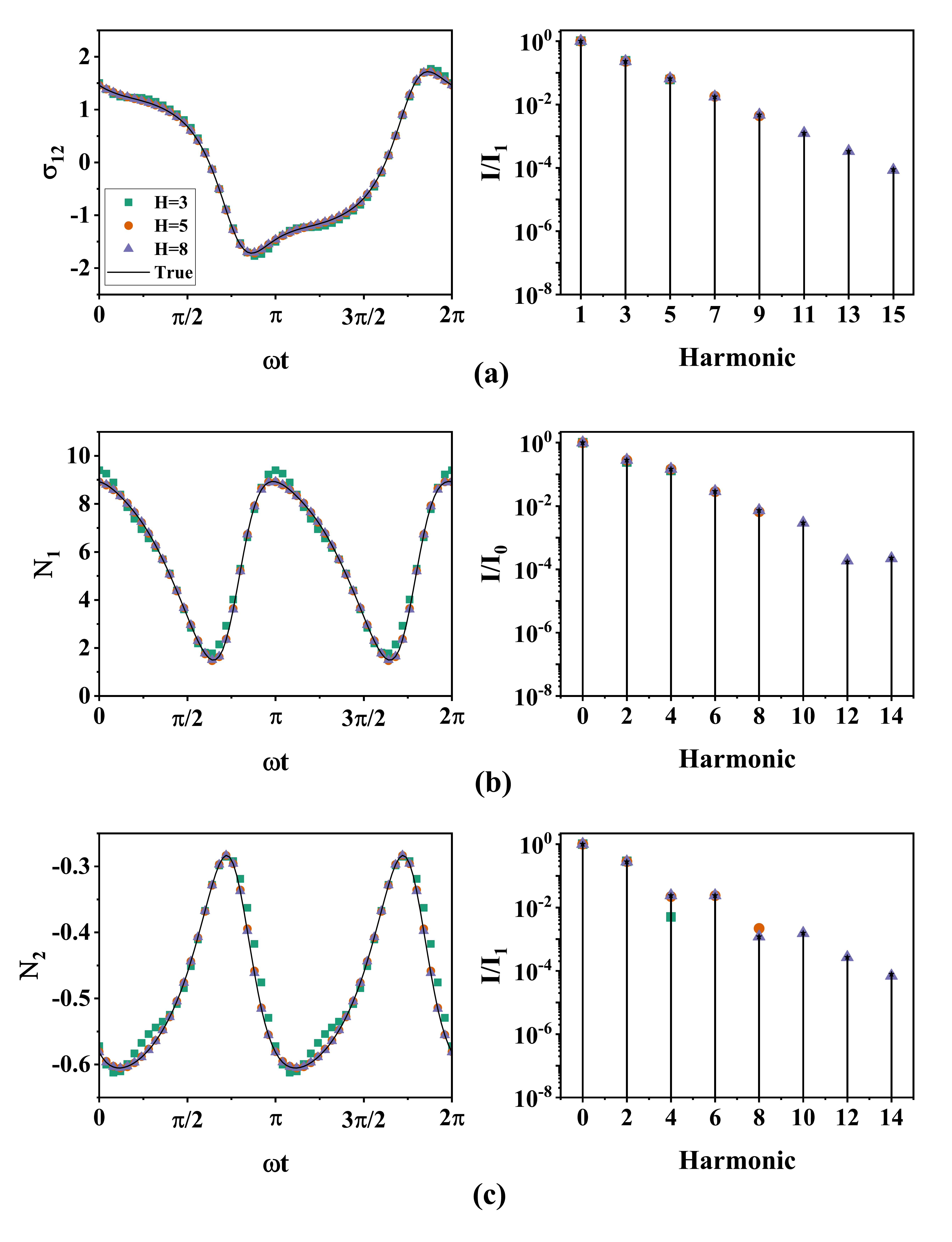}
\par\end{centering}
\caption{\label{fig:H_analysis_1}Stress waveforms and harmonic intensities
for (a) $\sigma_{12}$, (b) $N_{1}$ and (c) $N_{2}$ at different
$H$ for $\gamma_{0}=10$ and $\omega=1$ rad/s. The black solid lines
(left) and stars (right) correspond to the true solution, while the
symbols indicate results at different $H$. }
\end{figure}

\subsection{Convergence and the number of harmonics}

\begin{figure}
\begin{centering}
\includegraphics[scale=0.25]{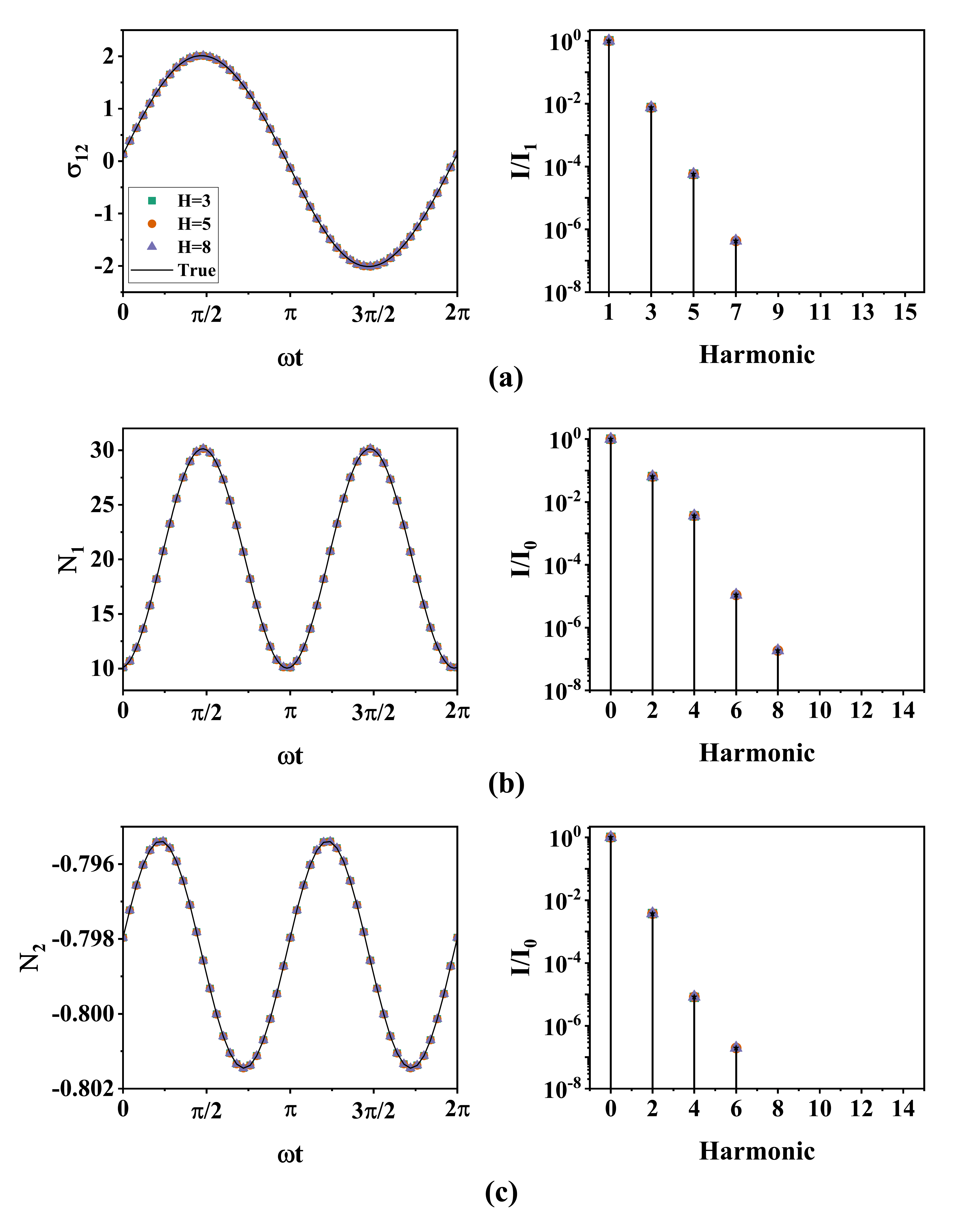}
\par\end{centering}
\caption{\label{fig:H_analysis_100}Stress waveforms and harmonic intensities
for (a) $\sigma_{12}$, (b) $N_{1}$ and (c) $N_{2}$ at different
$H$ for $\gamma_{0}=10$ and $\omega=100$ rad/s. The black solid
lines (left) and stars (right) correspond to the true solution, while
the symbols indicate results at different $H$.}
\end{figure}

A parsimonius choice for the number of harmonics $H$ ensures efficiency without compromising accuracy. Figures \ref{fig:H_analysis_1} and \ref{fig:H_analysis_100} show the stress waveforms and intensity plots for $\sigma_{12}$, $N_{1}$ and $N_{2}$ at $\gamma_{0}=10$ with $\omega=1$ rad/s and $\omega=100$ rad/s respectively, obtained using different values of $H$. The results are juxtaposed with the ``true" solution to the Giesekus model equations which is approximated using a large number of harmonics, viz. $H = 30$. From previous results, we know that higher harmonics decay slowly at large $\Wi$ and moderate $\De$. Therefore, selecting a large value for $H$ improves accuracy as shown in figure \ref{fig:H_analysis_1}. This impact is most pronounced in the second normal stress difference, figure \ref{fig:H_analysis_1}(c). In figure \ref{fig:H_analysis_100} in contrast, where $\gamma_0 = 10$ and $\De = 100$, higher harmonics decay rapidly, and increasing $H$ does not significantly improve accuracy.

Regardless, the magnitudes of the individual harmonic intensities are relatively stable as $H$ is increased. This provides a potential future pathway for systematically increasing the number of harmonics based on the decay characteristics of the harmonics, and level of accuracy desired. Such an adaptive protocol is not implemented in this study.

We now discuss the convergence of the HB method. Error bounds and convergence theorems for HB in different nonlinear vibration problems have already been established \cite{krack2019harmonic}. Here, we limit our discussion to harmonic convergence in systems that produce analytic outputs. Eqn. \eqref{eq:error}, which describes the pattern of convergence for such systems, can be adapted for our system using equations \eqref{eq:eq_for_T11}, \eqref{eq:eq_for_T22} and \eqref{eq:eq_for_T12} with $\bm{y}=[N_{1},  N_{2}, \sigma_{12}]^{T}$ and 
\begin{equation}
\begin{gathered}\xi{}_{F}=\left|\left|N_{1,H}-N_{1}(t)\right|\right|_{\infty}\leq u_{F}e^{-m_{F}\left(2H-2\right)},\\
\xi_{S}=\left|\left|N_{2,H}-N_{1}(t)\right|\right|_{\infty}\leq u_{S}e^{-m_{S}\left(2H-2\right),}\\
\xi_{G}=\left|\left|\sigma_{12,H}-\sigma_{12}(t)\right|\right|_{\infty}\leq u_{G}e^{-m_{G}\left(2H-1\right)},
\end{gathered}
\label{eq:norm}
\end{equation}
where $N_{1,H},N_{2,H}$ and $\sigma_{12,H}$ are HB approximations to the true solutions $N_1$, $N_2$, and $\sigma_{12}$ over a period of oscillation. The observed errors $\xi_{F},\xi_{S}$, and $\xi_{G}$ can be computed, and the decay coefficients $m_F$, $m_S$, and $m_G$ can be empirically studied as a function of $\alpha$, $\De$, and $\Wi$.

\begin{figure}
\begin{centering}
\includegraphics[scale=0.25]{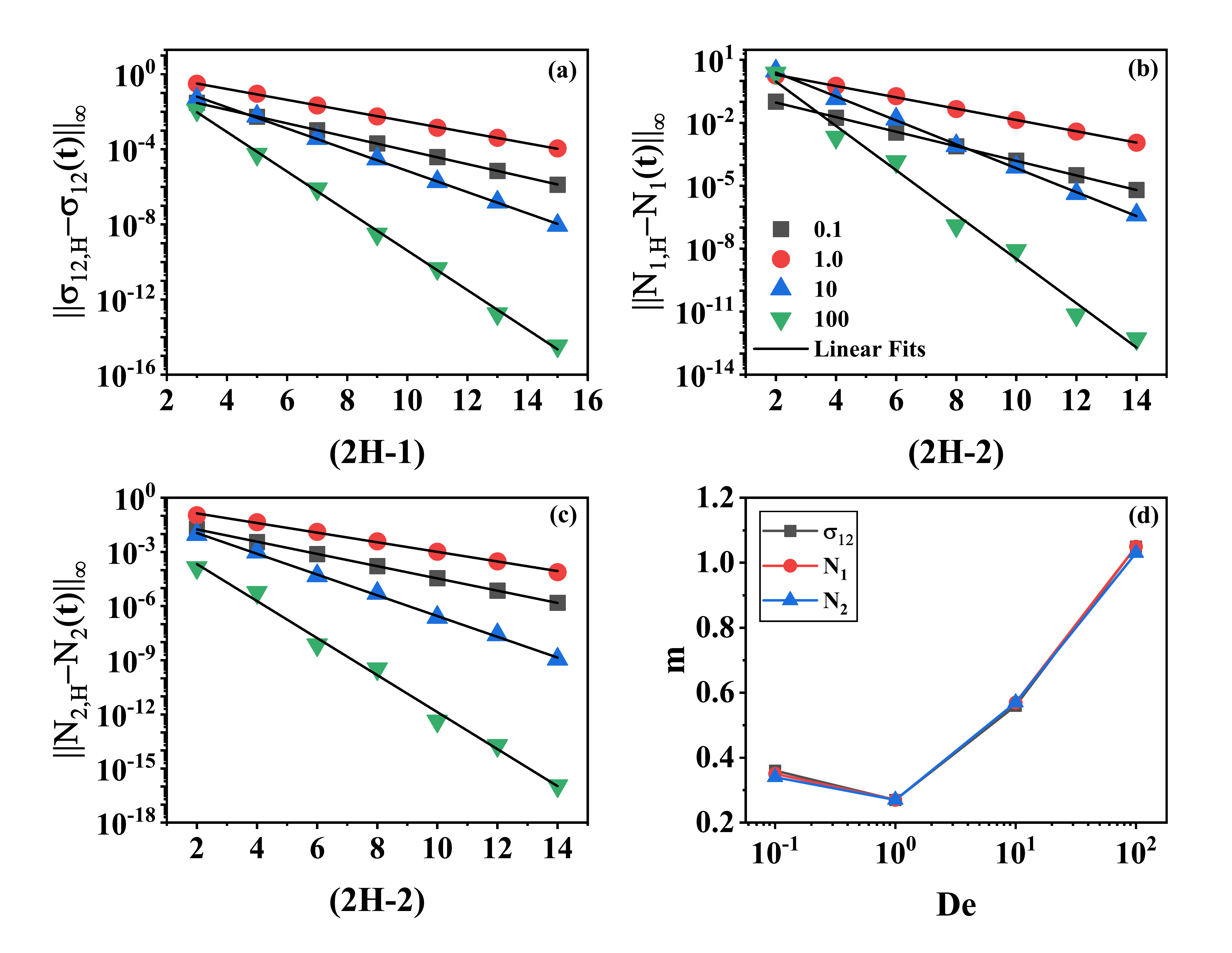}
\par\end{centering}
\caption{\label{fig:convergence} Convergence of the (a) shear stress,
(b) first normal stress difference, and (c) second normal
stress difference at different frequencies for $\gamma_{0}=10$, and
(d) the corresponding exponential decay coefficient, $m$. }
\end{figure}

Figure \ref{fig:convergence} shows $\xi_{F}$ and $\xi_{G}$ plotted against $(2H-1)$ and $(2H-2)$, respectively, at $\gamma_{0}=10$ for four different frequencies. The empirically obtained error in all the plots shows an exponential decay to the true solution. We observe that convergence is most rapid at high $\omega$, and most sluggish at intermediate $\omega$. Interestingly, the decay coefficient $m$ is almost the same for all three quantities i.e. $m_{F} \approx m_{S} \approx m_{G}$ as shown in figure \ref{fig:convergence}(d). 

\begin{figure}
\begin{centering}
\includegraphics[scale=0.3]{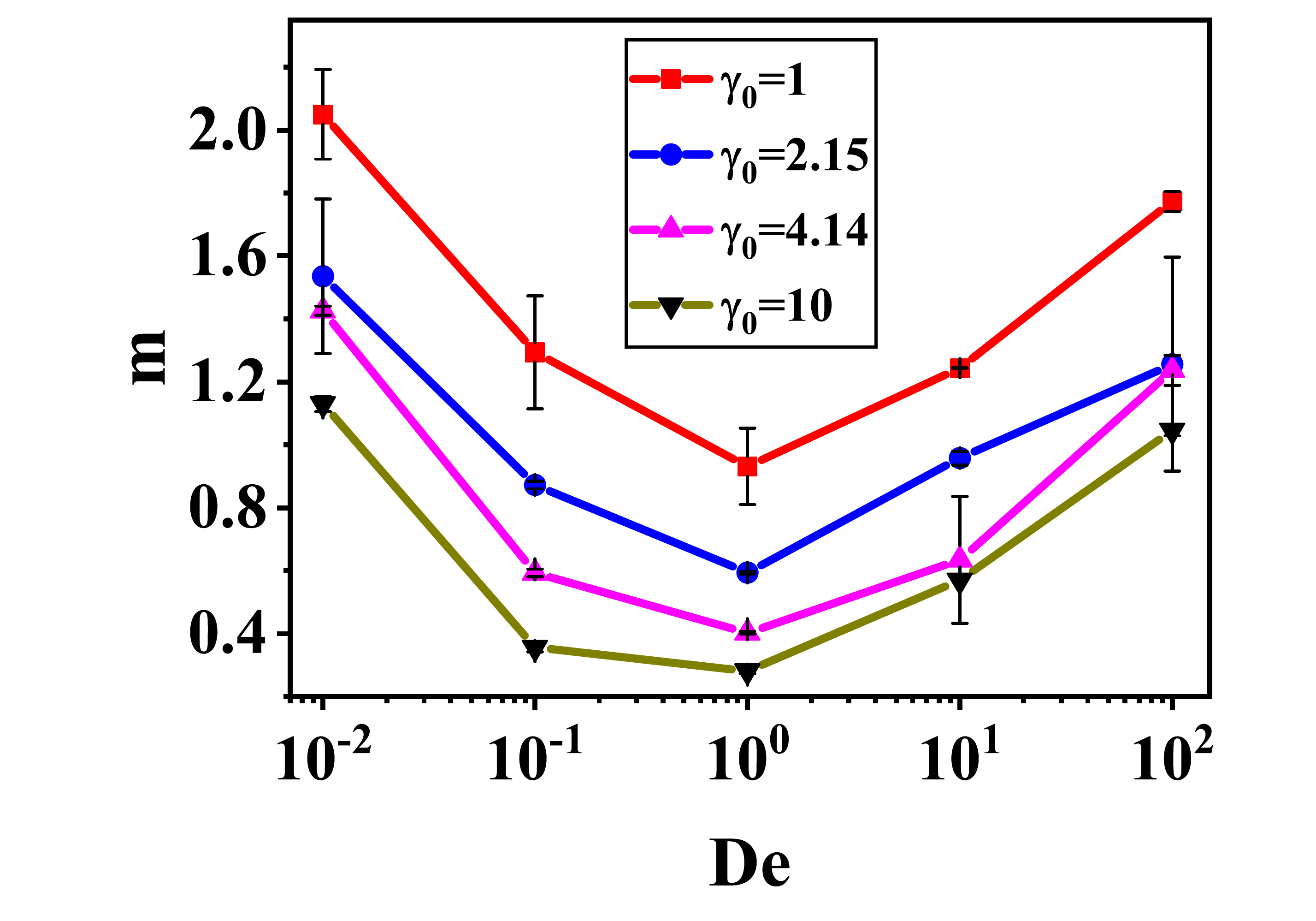}
\par\end{centering}
\caption{\label{fig:m_exponent}The convergence decay coefficient $m$ at different
strain amplitudes and frequencies.}
\end{figure}

We show the dependence of convergence rate on $\gamma_0$ in figure \ref{fig:m_exponent}. As expected, the pace of convergence slows down with increasing $\gamma_{0}$ at a particular frequency. The value of $m$ is inversely related to the extent of nonlinearity present in the system. The ``U" shape of these curves indicates the linearity of the system at the low and high frequency ends despite being deep into the LAOS regime at large $\gamma_{0}$. 

\subsection{Harmonic Balance vs Numerical Integration}

In this section, we demonstrate the compelling value proposition of HB. It is both faster and more accurate than the standard NI technique. The solid lines in figures \ref{fig:shear_coeffs}, \ref{fig:n1_coeffs}, and \ref{fig:n2_coeffs} show the numerically integrated solution to the IVP for the shear stress, first normal stress difference, and second normal stress difference respectively.
Visually, the agreement between the harmonic coefficients obtained
from both HB and numerical integration (NI) is excellent. Furthermore,
the agreement between these numerical solutions and the analytical
MAOS solution at $\gamma_{0}=0.1$ is also quite remarkable.

\begin{figure*}
\begin{center}
\includegraphics[scale=0.25]{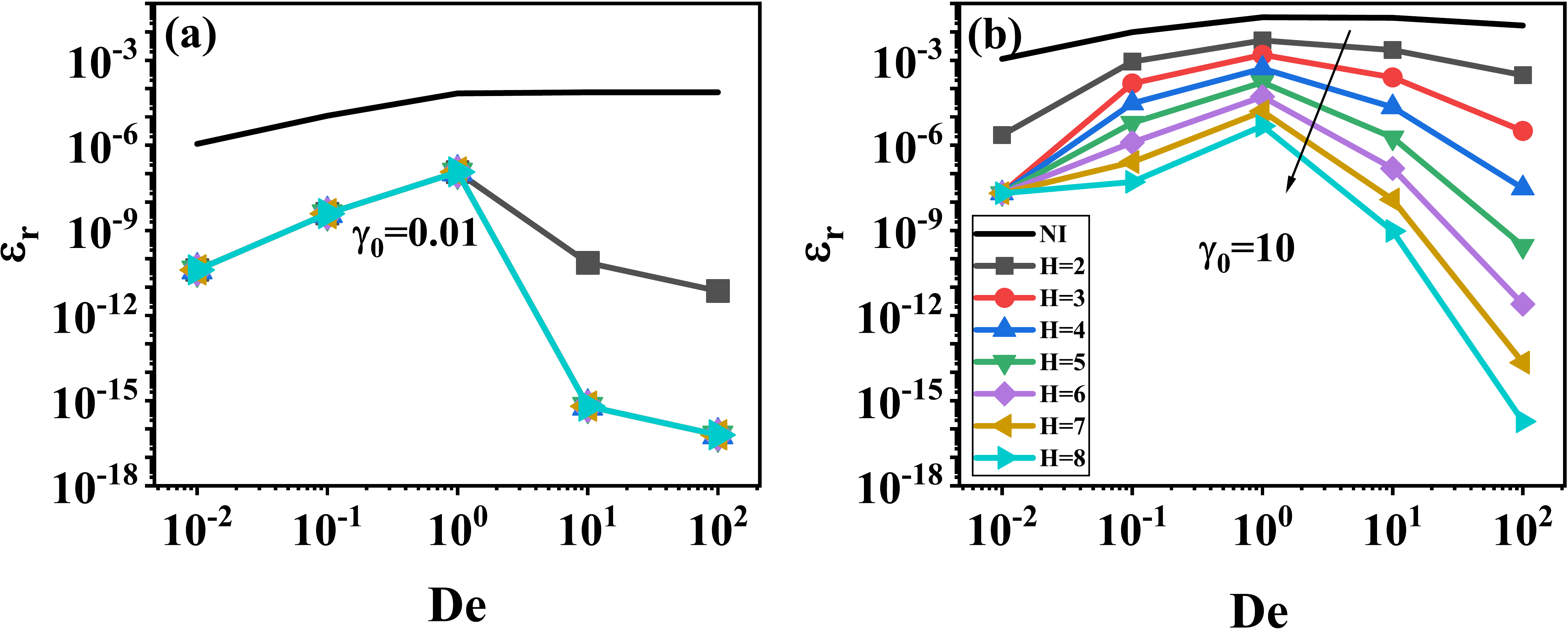}
\end{center}
\caption{\label{fig:rms_residual} The residual $\epsilon_r$ (eqn \ref{eq:residual2} ) for HB and NI at (a) $\gamma_{0}=0.01$ and (b) $\gamma_{0}=10$.}
\end{figure*}

Figure \ref{fig:rms_residual} depicts the residual $\epsilon_{r}$ (eqn \ref{eq:residual2}) at the two ends of the strain amplitude range, $\gamma_{0}=0.01$ and $\gamma_{0}=10$, for certain frequencies using both HB and NI. It is apparent from these plots that in terms of accuracy HB is superior to NI regardless of $\omega$ or $\gamma_0$. At $\gamma_{0} = 0.01$, only two harmonics ($H=2$) are required to achieve an $\epsilon_r$ that is at least three orders of magnitude smaller than NI. Using a single additional harmonic term ($H=3$) significantly improves the accuracy of HB at high frequencies at this $\gamma_{0}$. It is apparent from figure \ref{fig:rms_residual}(a) that increasing the number of harmonics any further does not improve the accuracy significantly. However, when $\gamma_0$ is increased to 10 in figure \ref{fig:rms_residual}(b),  increasing $H$ improves solution accuracy. However, even in this setting, HB with $H = 2$ harmonics outperforms NI in terms of accuracy.

\begin{figure*}
\begin{centering}
\includegraphics[scale=0.25]{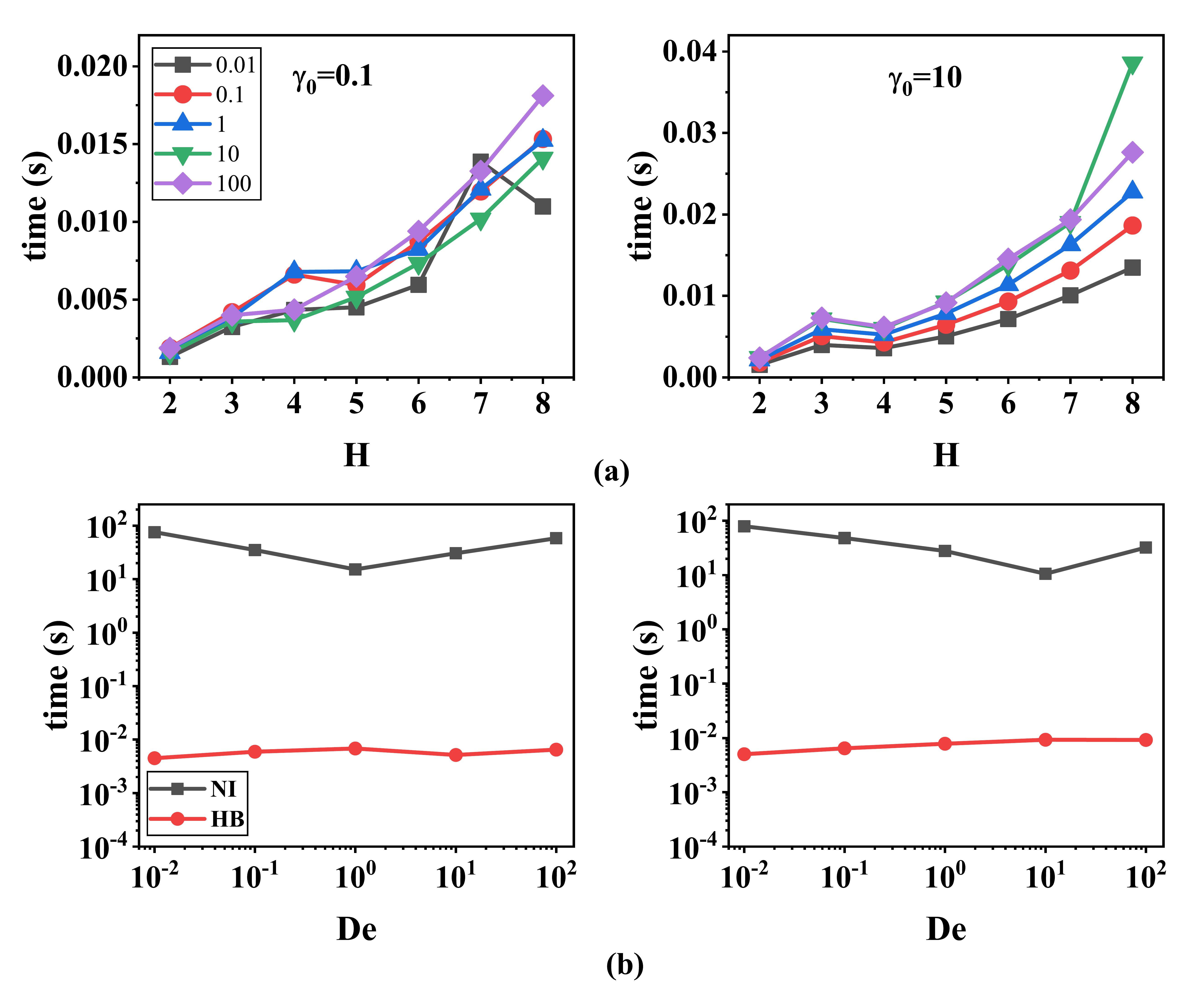} 
\par\end{centering}
\caption{\label{fig:cpu}(a) The computation time for the HB method at different
frequencies for $\gamma_{0}=0.1$ and $\gamma_{0}=10$. (b) Comparison
of CPU time for the numerically integrated IVP and HB with $H=5$
at different frequencies with $\gamma_{0}=0.1$ and $\gamma_{0}=10$.}
\end{figure*}

However it is not fair to only compare accuracy without acknowledging the computational cost incurred in producing more accurate solutions. HB generates a system of $(6H-2)$ nonlinear equations, which are evaluated at each iteration. Figure \ref{fig:cpu}(a) shows the computation time on a standard laptop computer (Intel Xeon W-1290 CPU 3.2GGHz) using HB for $\gamma_{0}=0.1$ and $\gamma_{0}=10$ at different frequencies. The number of equations, and the number of iterations required for convergence increase with $H$. The total computation time increases slightly faster than linearly with $H$. It also increases modestly with increasing $\omega$.

On the other hand, the accuracy and computational cost of NI is dominated by the computation of the periodic steady state solution. Figure \ref{fig:cpu}(b) compares the computational cost of NI and HB. NI takes roughly
10 -- 100s depending on the frequency. In contrast, the computational cost of HB is at least three orders of magnitude lower as shown in Figure \ref{fig:cpu}(b). Thus, HB is superior to NI in terms of accuracy, as well as speed.

\section{Conclusions and Perspective}

In this work, HB is used to solve the Giesekus model under oscillatory shear. HB can be thought of as a numerical extension to the analytical approach of finding higher-order Fourier coefficients for constitutive models subjected to LAOS. It directly computes the periodic steady state solution by reducing a system of ODEs to a system of nonlinear equations, and avoids solving an IVP in the time domain. Numerical approximations converge exponentially to the true solution as the number of harmonics probed is increased. This enables fast and accurate estimation of higher-order Fourier coefficients. Due to its speed and accuracy, HB facilitates parameter estimation and model selection using a Bayesian framework that enables interpretation of LAOS data through the lens of constitutive theories.

The nonlinearity in the Giesekus model is quadratic. HB handles such nonlinearities elegantly as illustrated in this work. However, HB is not limited to the Giesekus model, and can accommodate a wide variety of nonlinear models using the alternating frequency-time (AFT) technique \cite{krack2019harmonic}. Ongoing work suggests that AFT enables generalization of HB to arbitrary differential constitutive equations. In the future, numerical path continuation techniques can also be used to explore stable and unstable solution branches.


\section*{Acknowledgements}

This material is based partially upon work supported by the National Science Foundation under Grant No. DMR 1727870 (SS). Financial support from the Science and Engineering Research Board, Government of India and Prime Minister Research Fellowship, Ministry of Education, Government of India is also acknowledged (YMJ and SM). 



\bibliographystyle{elsarticle-num} 
\bibliography{HB_references}

\begin{thebibliography}{10}
\expandafter\ifx\csname url\endcsname\relax
  \def\url#1{\texttt{#1}}\fi
\expandafter\ifx\csname urlprefix\endcsname\relax\def\urlprefix{URL }\fi
\expandafter\ifx\csname href\endcsname\relax
  \def\href#1#2{#2} \def\path#1{#1}\fi

\bibitem{ma2020effects}
Y.~Ma, D.~Su, Y.~Wang, D.~Li, L.~Wang, Effects of concentration and nacl on
  rheological behaviors of konjac glucomannan solution under large amplitude
  oscillatory shear ({LAOS}), LWT 128 (2020) 109466.

\bibitem{khandavalli2015large}
S.~Khandavalli, J.~P. Rothstein, Large amplitude oscillatory shear rheology of
  three different shear-thickening particle dispersions, Rheol. Acta 54~(7)
  (2015) 601--618.

\bibitem{chan2018nonlinear}
R.~W. Chan, Nonlinear viscoelastic characterization of human vocal fold tissues
  under large-amplitude oscillatory shear ({LAOS}), J. Rheol. 62~(3) (2018)
  695--712.

\bibitem{wapperom2005numerical}
P.~Wapperom, A.~Leygue, R.~Keunings, P.~Wapperom, A.~Leygue, R.~Keunings,
  Numerical simulation of large amplitude oscillatory shear of a high-density
  polyethylene melt using the {MSF} model, J. Non-Newtonian Fluid Mech. 130~(2)
  (2005) 63--76.
\newblock \href {https://doi.org/10.1016/j.jnnfm.2005.08.002}
  {\path{doi:10.1016/j.jnnfm.2005.08.002}}.

\bibitem{li2009nonlinearity}
X.~Li, S.-Q. Wang, X.~Wang, Nonlinearity in large amplitude oscillatory shear
  ({LAOS}) of different viscoelastic materials, J. Rheol. 53~(5) (2009)
  1255--1274.

\bibitem{dimitriou2013describing}
C.~J. Dimitriou, R.~H. Ewoldt, G.~H. McKinley, Describing and prescribing the
  constitutive response of yield stress fluids using large amplitude
  oscillatory shear stress ({LAOStress}), J. Rheol. 57~(1) (2013) 27--70.

\bibitem{min2014microstructure}
J.~Min~Kim, A.~P. Eberle, A.~Kate~Gurnon, L.~Porcar, N.~J. Wagner, The
  microstructure and rheology of a model, thixotropic nanoparticle gel under
  steady shear and large amplitude oscillatory shear ({LAOS}), J. Rheol. 58~(5)
  (2014) 1301--1328.

\bibitem{armstrong2016dynamic}
M.~J. Armstrong, A.~N. Beris, S.~A. Rogers, N.~J. Wagner, Dynamic shear
  rheology of a thixotropic suspension: Comparison of an improved
  structure-based model with large amplitude oscillatory shear experiments, J.
  Rheol. 60~(3) (2016) 433--450.

\bibitem{armstrong2021simple}
M.~Armstrong, M.~Scully, M.~Clark, T.~Corrigan, C.~James, A simple approach for
  adding thixotropy to an elasto-visco-plastic rheological model to facilitate
  structural interrogation of human blood, J. Non-Newtonian Fluid Mech. 290
  (2021) 104503.
\newblock \href {https://doi.org/10.1016/j.jnnfm.2021.104503}
  {\path{doi:10.1016/j.jnnfm.2021.104503}}.

\bibitem{Donley2019}
G.~J. Donley, J.~R. {de Bruyn}, G.~H. McKinley, S.~A. Rogers, Time-resolved
  dynamics of the yielding transition in soft materials, J. Non-Newtonian Fluid
  Mech. 264 (2019) 117--134.
\newblock \href {https://doi.org/10.1016/j.jnnfm.2018.10.003}
  {\path{doi:10.1016/j.jnnfm.2018.10.003}}.

\bibitem{ewoldt2010large}
R.~H. Ewoldt, P.~Winter, J.~Maxey, G.~H. McKinley, Large amplitude oscillatory
  shear of pseudoplastic and elastoviscoplastic materials, Rheol. Acta 49~(2)
  (2010) 191--212.

\bibitem{stickel2013response}
J.~J. Stickel, J.~S. Knutsen, M.~W. Liberatore, Response of elastoviscoplastic
  materials to large amplitude oscillatory shear flow in the parallel-plate and
  cylindrical-couette geometries, J. Rheol. 57~(6) (2013) 1569--1596.

\bibitem{dimitriou2012rheo}
C.~J. Dimitriou, L.~Casanellas, T.~J. Ober, G.~H. McKinley, Rheo-piv of a
  shear-banding wormlike micellar solution under large amplitude oscillatory
  shear, Rheol. Acta 51~(5) (2012) 395--411.

\bibitem{goudoulas2017nonlinearities}
T.~B. Goudoulas, S.~Pan, N.~Germann, Nonlinearities and shear banding
  instability of polyacrylamide solutions under large amplitude oscillatory
  shear, J. Rheol. 61~(5) (2017) 1061--1083.

\bibitem{kate2012large}
A.~Kate~Gurnon, N.~J. Wagner, Large amplitude oscillatory shear ({LAOS})
  measurements to obtain constitutive equation model parameters: {Giesekus}
  model of banding and nonbanding wormlike micelles, J. Rheol. 56~(2) (2012)
  333--351.

\bibitem{radhakrishnan2018shear}
R.~Radhakrishnan, S.~M. Fielding, Shear banding in large amplitude oscillatory
  shear ({LAOStrain} and {LAOStress}) of soft glassy materials, J. Rheol.
  62~(2) (2018) 559--576.

\bibitem{atalik2004occurrence}
K.~Atal{\i}k, R.~Keunings, On the occurrence of even harmonics in the shear
  stress response of viscoelastic fluids in large amplitude oscillatory shear,
  J. Non-Newtonian Fluid Mech. 122~(1-3) (2004) 107--116.

\bibitem{yang2017dynamic}
K.~Yang, W.~Yu, Dynamic wall slip behavior of yield stress fluids under large
  amplitude oscillatory shear., J. Rheol. 61~(4) (2017) 627--641.

\bibitem{klein2007separation}
C.~O. Klein, H.~W. Spiess, A.~Calin, C.~Balan, M.~Wilhelm, Separation of the
  nonlinear oscillatory response into a superposition of linear, strain
  hardening, strain softening, and wall slip response, Macromolecules 40~(12)
  (2007) 4250--4259.

\bibitem{graham1995wall}
M.~D. Graham, Wall slip and the nonlinear dynamics of large amplitude
  oscillatory shear flows, J. Rheol. 39~(4) (1995) 697--712.

\bibitem{suman2022large}
K.~Suman, S.~Shanbhag, Y.~M. Joshi, Large amplitude oscillatory shear study of
  a colloidal gel at the critical state, arXiv preprint arXiv:2211.16724
  (2022).

\bibitem{kim2014microstructure}
J.~Kim, D.~Merger, M.~Wilhelm, M.~E. Helgeson, Microstructure and nonlinear
  signatures of yielding in a heterogeneous colloidal gel under large amplitude
  oscillatory shear, J. Rheol. 58~(5) (2014) 1359--1390.

\bibitem{sun2015large}
W.-x. Sun, L.-z. Huang, Y.-r. Yang, X.-x. Liu, Z.~Tong, Large amplitude
  oscillatory shear studies on the strain-stiffening behavior of gelatin gels,
  Chin. J. Polym. Sci. 33~(1) (2015) 70--83.

\bibitem{ng2011large}
T.~S. Ng, G.~H. McKinley, R.~H. Ewoldt, Large amplitude oscillatory shear flow
  of gluten dough: A model power-law gel, J. Rheol. 55~(3) (2011) 627--654.

\bibitem{wagner2011analysis}
M.~H. Wagner, V.~H. Rol{\'o}n-Garrido, K.~Hyun, M.~Wilhelm, Analysis of medium
  amplitude oscillatory shear data of entangled linear and model comb polymers,
  J. Rheol. 55~(3) (2011) 495--516.

\bibitem{hyun2007fourier}
K.~Hyun, E.~S. Baik, K.~H. Ahn, S.~J. Lee, M.~Sugimoto, K.~Koyama,
  Fourier-transform rheology under medium amplitude oscillatory shear for
  linear and branched polymer melts, J. Rheol. 51~(6) (2007) 1319--1342.

\bibitem{hoyle2014large}
D.~Hoyle, D.~Auhl, O.~Harlen, V.~Barroso, M.~Wilhelm, T.~McLeish, Large
  amplitude oscillatory shear and fourier transform rheology analysis of
  branched polymer melts, J. Rheol. 58~(4) (2014) 969--997.

\bibitem{cho2010scaling}
K.~S. Cho, K.-W. Song, G.-S. Chang, Scaling relations in nonlinear viscoelastic
  behavior of aqueous peo solutions under large amplitude oscillatory shear
  flow, Journal of Rheology 54~(1) (2010) 27--63.
\newblock \href {https://doi.org/10.1122/1.3258278}
  {\path{doi:10.1122/1.3258278}}.

\bibitem{cho2015effect}
K.~S. Cho, J.~W. Kim, J.-E. Bae, J.~H. Youk, H.~J. Jeon, K.-W. Song, Effect of
  temporary network structure on linear and nonlinear viscoelasticity of
  polymer solutions, Korea-Aust. Rheol. J. 27~(2) (2015) 151--161.
\newblock \href {https://doi.org/10.1007/s13367-015-0015-y}
  {\path{doi:10.1007/s13367-015-0015-y}}.

\bibitem{cho2016viscoelasticity}
K.~S. Cho, Viscoelasticity of Polymers: {T}heory and Numerical Algorithms,
  Springer, Dordrecht, the Netherlands, 2016.

\bibitem{macdonald1969rheological}
I.~F. MacDonald, B.~D. Marsh, E.~Ashare, Rheological behavior for large
  amplitude oscillatory motion, Chem. Eng. Sci. 24~(10) (1969) 1615--1625.

\bibitem{pearson1982behavior}
D.~S. Pearson, W.~E. Rochefort, Behavior of concentrated polystyrene solutions
  in large-amplitude oscillating shear fields, J. Polym. Sci. Polym. Phys.
  20~(1) (1982) 83--98.

\bibitem{giacomin2015pade}
A.~J. Giacomin, C.~Saengow, M.~Guay, C.~Kolitawong, Pad{\'e} approximants for
  large-amplitude oscillatory shear flow, Rheol. Acta 54~(8) (2015) 679--693.

\bibitem{Ewoldt2008}
R.~H. Ewoldt, A.~Hosoi, G.~H. McKinley, New measures for characterizing
  nonlinear viscoelasticity in large amplitude oscillatory shear, J. Rheol.
  52~(6) (2008) 1427--1458.

\bibitem{cho2005geometrical}
K.~S. Cho, K.~Hyun, K.~H. Ahn, S.~J. Lee, A geometrical interpretation of large
  amplitude oscillatory shear response, J. Rheol. 49~(3) (2005) 747--758.

\bibitem{bae2017analytical}
J.-E. Bae, K.~S. Cho, Analytical studies on the {LAOS} behaviors of some
  popularly used viscoelastic constitutive equations with a new insight on
  stress decomposition of normal stresses, Phys. Fluids 29~(9) (2017) 093103.
\newblock \href {https://doi.org/10.1063/1.5001742}
  {\path{doi:10.1063/1.5001742}}.

\bibitem{hyun2009establishing}
K.~Hyun, M.~Wilhelm, Establishing a new mechanical nonlinear coefficient q from
  ft-rheology: First investigation of entangled linear and comb polymer model
  systems, Macromolecules 42~(1) (2009) 411--422.

\bibitem{rogers2012sequence}
S.~A. Rogers, M.~P. Lettinga, A sequence of physical processes determined and
  quantified in large-amplitude oscillatory shear ({LAOS}): {Application} to
  theoretical nonlinear models, J Rheol 56~(1) (2012) 1--25.

\bibitem{ferry1980viscoelastic}
J.~D. Ferry, Viscoelastic properties of polymers, $3^\text{rd}$ Edition, John
  Wiley \& Sons, New York, NY, 1980.

\bibitem{giacomin1993large}
A.~J. Giacomin, J.~M. Dealy, Large-amplitude oscillatory shear, in: A.~A.
  Collyer (Ed.), Techniques in Rheological Measurement, Springer Netherlands,
  Dordrecht, 1993, pp. 99--121.

\bibitem{nam2008prediction}
J.~G. Nam, K.~Hyun, K.~H. Ahn, S.~J. Lee, Prediction of normal stresses under
  large amplitude oscillatory shear flow, J. Non-Newtonian Fluid Mech. 150~(1)
  (2008) 1--10.

\bibitem{shanbhag2022kramers1}
S.~Shanbhag, Y.~M. Joshi, Kramers--kronig relations for nonlinear rheology.
  part i: General expression and implications, J. Rheol. 66~(5) (2022)
  973--982.

\bibitem{shanbhag2022kramers2}
S.~Shanbhag, Y.~M. Joshi, Kramers--kronig relations for nonlinear rheology.
  part ii: Validation of medium amplitude oscillatory shear (maos)
  measurements, J. Rheol. 66~(5) (2022) 925--936.

\bibitem{saengow2015exact}
C.~Saengow, A.~J. Giacomin, C.~Kolitawong, Exact analytical solution for
  large-amplitude oscillatory shear flow, Macromol. Theory Simul. 24~(4) (2015)
  352--392.

\bibitem{poungthong2019exact}
P.~Poungthong, A.~Giacomin, C.~Saengow, C.~Kolitawong, Exact solution for
  intrinsic nonlinearity in oscillatory shear from the corotational maxwell
  fluid, J. Non-Newtonian Fluid Mech. 265 (2019) 53--65.
\newblock \href {https://doi.org/10.1016/j.jnnfm.2019.01.001}
  {\path{doi:10.1016/j.jnnfm.2019.01.001}}.

\bibitem{saengow2017exact}
C.~Saengow, A.~J. Giacomin, C.~Kolitawong, Exact analytical solution for
  large-amplitude oscillatory shear flow from oldroyd 8-constant framework:
  Shear stress, Phys. Fluids 29~(4) (2017) 043101.

\bibitem{saengow2017normal}
C.~Saengow, A.~J. Giacomin, Normal stress differences from oldroyd 8-constant
  framework: Exact analytical solution for large-amplitude oscillatory shear
  flow, Phys. Fluids 29~(12) (2017) 121601.

\bibitem{helfand1982calculation}
E.~Helfand, D.~S. Pearson, Calculation of the nonlinear stress of polymers in
  oscillatory shear fields, J. Polym. Sci. Polym. Phys. 20~(7) (1982)
  1249--1258.

\bibitem{bharadwaj2015constitutive}
N.~A. Bharadwaj, R.~H. Ewoldt, Constitutive model fingerprints in
  medium-amplitude oscillatory shear, J. Rheol. 59~(2) (2015) 557--592.

\bibitem{bird2014dilute}
R.~Bird, A.~Giacomin, A.~Schmalzer, C.~Aumnate, Dilute rigid dumbbell
  suspensions in large-amplitude oscillatory shear flow: Shear stress response,
  J. Chem. Phys. 140~(7) (2014) 074904.

\bibitem{fan1984kinetic}
X.-J. Fan, R.~B. Bird, A kinetic theory for polymer melts vi. calculation of
  additional material functions, J. Non-Newtonian Fluid Mech. 15~(3) (1984)
  341--373.

\bibitem{bharadwaj2014general}
N.~A. Bharadwaj, R.~H. Ewoldt, The general low-frequency prediction for
  asymptotically nonlinear material functions in oscillatory shear, J. Rheol.
  58~(4) (2014) 891--910.

\bibitem{yu2002modeling}
W.~Yu, M.~Bousmina, M.~Grmela, C.~Zhou, Modeling of oscillatory shear flow of
  emulsions under small and large deformation fields, J. Rheol. 46~(6) (2002)
  1401--1418.

\bibitem{martinetti2019time}
L.~Martinetti, R.~H. Ewoldt, Time-strain separability in medium-amplitude
  oscillatory shear, Phys. Fluids 31~(2) (2019) 021213.

\bibitem{shanbhag2021spectral}
S.~Shanbhag, S.~Mittal, Y.~M. Joshi, Spectral method for time-strain separable
  integral constitutive models in oscillatory shear, Phys. Fluids 33~(11)
  (2021) 113104.

\bibitem{Shanbhag2010}
S.~Shanbhag, Analytical rheology of blends of linear and star polymers using a
  {Bayesian} formulation, Rheol. Acta 49~(4) (2010) 411--422.

\bibitem{Takeh2011}
A.~Takeh, J.~Worch, S.~Shanbhag, Analytical rheology of metallocene-catalyzed
  polyethylenes, Macromolecules 44~(9) (2011) 3656--3665.
\newblock \href {https://doi.org/10.1021/ma2004772}
  {\path{doi:10.1021/ma2004772}}.

\bibitem{giesekus1982simple}
H.~Giesekus, A simple constitutive equation for polymer fluids based on the
  concept of deformation-dependent tensorial mobility, J. Non-Newtonian Fluid
  Mech. 11~(1-2) (1982) 69--109.

\bibitem{yoo1989steady}
J.~Yoo, H.~C. Choi, On the steady simple shear flows of the one-mode {Giesekus}
  fluid, Rheol. Acta 28~(1) (1989) 13--24.

\bibitem{schleiniger1991remark}
G.~Schleiniger, R.~Weinacht, A remark on the {Giesekus} viscoelastic fluid, J.
  Rheol. 35~(6) (1991) 1157--1170.

\bibitem{yao1998extensional}
M.~Yao, G.~H. McKinley, B.~Debbaut, Extensional deformation, stress relaxation
  and necking failure of viscoelastic filaments, J. Non-Newtonian Fluid Mech.
  79~(2-3) (1998) 469--501.

\bibitem{holz1999shear}
T.~Holz, P.~Fischer, H.~Rehage, Shear relaxation in the nonlinear-viscoelastic
  regime of a {Giesekus} fluid, J. Non-Newtonian Fluid Mech. 88~(1-2) (1999)
  133--148.

\bibitem{fischer1997non}
P.~Fischer, H.~Rehage, Non-linear flow properties of viscoelastic surfactant
  solutions, Rheol. Acta 36~(1) (1997) 13--27.

\bibitem{rehage2015experimental}
H.~Rehage, R.~Fuchs, Experimental and numerical investigations of the
  non-linear rheological properties of viscoelastic surfactant solutions:
  application and failing of the one-mode {Giesekus} model, Colloid Polym. Sci.
  293~(11) (2015) 3249--3265.

\bibitem{bandyopadhyay2005effect}
R.~Bandyopadhyay, A.~Sood, Effect of silica colloids on the rheology of
  viscoelastic gels formed by the surfactant cetyl trimethylammonium tosylate,
  J. Colloid Interface Sci. 283~(2) (2005) 585--591.

\bibitem{kokini2000integral}
J.~Kokini, M.~Dhanasekharan, C.~Wang, H.~Huang, Integral and differential
  linear and non-linear constitutive models for rheology of wheat flour doughs,
  Technomics Publishing Co. Inc., Lancaster, PA, 2000.

\bibitem{Dhanasekharan2001}
M.~Dhanasekharan, C.~Wang, J.~Kokini, Use of nonlinear differential
  viscoelastic models to predict the rheological properties of gluten dough, J.
  Food Process Eng 24~(3) (2001) 193--216.

\bibitem{duvarci2017saos}
O.~C. Duvarci, G.~Yazar, J.~L. Kokini, The {SAOS}, {MAOS} and {LAOS} behavior
  of a concentrated suspension of tomato paste and its prediction using the
  {Bird-Carreau} ({SAOS}) and {Giesekus} models ({MAOS-LAOS}), J. Food Eng. 208
  (2017) 77--88.

\bibitem{calin2010determination}
A.~Calin, M.~Wilhelm, C.~Balan, Determination of the non-linear parameter
  (mobility factor) of the {Giesekus} constitutive model using {LAOS}
  procedure, J. Non-Newtonian Fluid Mech. 165~(23-24) (2010) 1564--1577.

\bibitem{quinzani1990modeling}
L.~Quinzani, G.~McKinley, R.~Brown, R.~Armstrong, Modeling the rheology of
  polyisobutylene solutions, J. Rheol. 34~(5) (1990) 705--748.

\bibitem{debbaut2002large}
B.~Debbaut, H.~Burhin, Large amplitude oscillatory shear and fourier-transform
  rheology for a high-density polyethylene: Experiments and numerical
  simulation, J. Rheol. 46~(5) (2002) 1155--1176.

\bibitem{oztekin1994quantitative}
A.~{\"O}ztekin, R.~A. Brown, G.~H. McKinley, Quantitative prediction of the
  viscoelastic instability in cone-and-plate flow of a {Boger} fluid using a
  multi-mode {Giesekus} model, J. Non-Newtonian Fluid Mech. 54 (1994) 351--377.

\bibitem{atalik2002non}
K.~Atal{\i}k, R.~Keunings, Non-linear temporal stability analysis of
  viscoelastic plane channel flows using a fully-spectral method, J.
  Non-Newtonian Fluid Mech. 102~(2) (2002) 299--319.
\newblock \href {https://doi.org/10.1016/S0377-0257(01)00184-7}
  {\path{doi:10.1016/S0377-0257(01)00184-7}}.

\bibitem{borzacchiello2016orientation}
D.~Borzacchiello, E.~Abisset-Chavanne, F.~Chinesta, R.~Keunings, Orientation
  kinematics of short fibres in a second-order viscoelastic fluid, Rheol. Acta
  55~(5) (2016) 397--409.
\newblock \href {https://doi.org/10.1007/s00397-016-0929-4}
  {\path{doi:10.1007/s00397-016-0929-4}}.

\bibitem{ewoldt2010secondary}
R.~H. Ewoldt, G.~H. McKinley, On secondary loops in {LAOS} via
  self-intersection of lissajous--bowditch curves, Rheol. Acta 49~(2) (2010)
  213--219.

\bibitem{bae2015semianalytical}
J.-E. Bae, K.~S. Cho, Semianalytical methods for the determination of the
  nonlinear parameter of nonlinear viscoelastic constitutive equations from
  {LAOS} data, J. Rheol. 59~(2) (2015) 525--555.

\bibitem{Atkinson2009}
K.~Atkinson, W.~Han, D.~E. Stewart, Numerical Solution of Ordinary Differential
  Equations, John Wiley \& Sons, Ltd, Hoboken, New Jersey, 2009.
\newblock \href {https://doi.org/https://doi.org/10.1002/9781118164495}
  {\path{doi:https://doi.org/10.1002/9781118164495}}.

\bibitem{krack2019harmonic}
M.~Krack, J.~Gross, Harmonic Balance for Nonlinear Vibration Problems, Springer
  International Publishing, Cham, 2019.
\newblock \href {https://doi.org/10.1007/978-3-030-14023-6}
  {\path{doi:10.1007/978-3-030-14023-6}}.

\bibitem{khair2016large}
A.~S. Khair, Large amplitude oscillatory shear of the giesekus model, J. Rheol.
  60~(2) (2016) 257--266.
\newblock \href {https://doi.org/10.1122/1.4941423}
  {\path{doi:10.1122/1.4941423}}.

\bibitem{jeyaseelan2008network}
R.~S. Jeyaseelan, A.~J. Giacomin, Network theory for polymer solutions in large
  amplitude oscillatory shear, J. Non-Newtonian Fluid Mech. 148~(1-3) (2008)
  24--32.

\bibitem{tee1975nonlinear}
T.-T. Tee, J.~Dealy, Nonlinear viscoelasticity of polymer melts, Trans. Soc.
  Rheol. 19~(4) (1975) 595--615.

\bibitem{khandavalli2016comparison}
S.~Khandavalli, J.~Hendricks, C.~Clasen, J.~P. Rothstein, A comparison of
  linear and branched wormlike micelles using large amplitude oscillatory shear
  and orthogonal superposition rheology, J. Rheol. 60~(6) (2016) 1331--1346.

\bibitem{moud2021viscoelastic}
A.~A. Moud, M.~Kamkar, A.~Sanati-Nezhad, S.~H. Hejazi, U.~Sundararaj,
  Viscoelastic properties of poly (vinyl alcohol) hydrogels with cellulose
  nanocrystals fabricated through sodium chloride addition: Rheological
  evidence of double network formation, Colloids Surf., A 609 (2021) 125577.

\bibitem{yazar2017non}
G.~Yazar, O.~Duvarci, S.~Tavman, J.~L. Kokini, Non-linear rheological behavior
  of gluten-free flour doughs and correlations of {LAOS} parameters with
  gluten-free bread properties, J. Cereal Sci. 74 (2017) 28--36.

\bibitem{yasin2021large}
S.~Yasin, M.~Hussain, Q.~Zheng, Y.~Song, Large amplitude oscillatory rheology
  of silica and cellulose nanocrystals filled natural rubber compounds, J.
  Colloid Interface Sci. 588 (2021) 602--610.

\bibitem{leygue2006tube}
A.~Leygue, C.~Bailly, R.~Keunings, A tube-based constitutive equation for
  polydisperse entangled linear polymers, J. Non-Newtonian Fluid Mech. 136~(1)
  (2006) 1--16.

\bibitem{hyun2013numerical}
K.~Hyun, W.~Kim, S.~Joon~Park, M.~Wilhelm, Numerical simulation results of the
  nonlinear coefficient {Q} from {FT}-rheology using a single mode pom-pom
  model, J. Rheol. 57~(1) (2013) 1--25.

\end{thebibliography}


%
%
%

\end{document}